\documentclass[11pt,a4paper]{amsart}
\usepackage[english]{babel}
\usepackage{amssymb, latexsym, amsthm,amsmath,tikz,verbatim,epigraph,graphicx} 
\usepackage{geometry}
\usepackage{enumitem}
\usepackage{tikz}
\usetikzlibrary{positioning}
\usetikzlibrary{matrix}
\usetikzlibrary{graphs}
\usetikzlibrary{trees}

\newtheorem{thm}{Theorem}[section]

\newtheorem{cor}[thm]{Corollary}
\newtheorem{defn}[thm]{Definition}

\newtheorem{exmpl}[thm]{Example}
\newtheorem{lem}[thm]{Lemma}
\newtheorem{prop}[thm]{Proposition}

\newtheorem{rem}[thm]{Remark}

\newtheorem{conc}[thm]{Conclusion}

\def\fpbar{\overline{\mathbb F}_p} 
\def\Fq{{\mathbb F}_q}
\def\indkzg{\mathrm{ind}^G_{KZ}} %

\def\gn{g^0_{n,\mu}}
\def\gnn{g^1_{n,\mu}}
\def\xvec{\bigotimes\limits_{j=0}^{f-1} x_j^{r_j}}
\def\xbutkvec{\bigotimes\limits_{k\neq j=0}^{f-1} x_j^{r_j}}
\def\ybutkvec{\bigotimes\limits_{k\neq j=0}^{f-1} y_j^{r_j}}
\def\xbutsvec{\bigotimes\limits_{s\neq j=0}^{f-1} x_j^{r_j}}

\def\yvec{\bigotimes\limits_{j=0}^{f-1} y_j^{r_j}}
\def \xyvec{\bigotimes \limits _{j=0}^{f-1} x_j^{r_j-i_j}y_j^{i_j}}
\def \yxvec{\bigotimes \limits _{j=0}^{f-1} x_j^{i_j}y_j^{r_j-i_j}}
\def\serwgt{\det^\omega \otimes\bigotimes\limits_{j=0}^{f-1}\mathrm{Sym}_j^{r_j}\fpbar^2}
\def\serwgtnew{\det^\omega \otimes\mathrm{Sym}^{\vec r}\fpbar^2}
\def\An{A_n^0(\sigma)}
\def\Ann{A_n^1(\sigma)}
\def\Qp{\mathbb Q_p}

\def\Ind{\mathrm {Ind}}
\def\Id{\mathrm {Id}}
\def\qBn{\overline{B}_n}

\newif\ifXY 
\XYtrue     
\ifXY
\usepackage{xy}
\xyoption{all} \fi

\AtEndDocument{\bigskip{\footnotesize%
  \textsc{Yotam Hendel, Faculty of Mathematics and Computer Science, Weizmann Institute of
Science,  Rehovot 76100, Israel.
} \par  
  \textit{E-mail address}: \texttt{yotam.hendel@gmail.com} \par
  \addvspace{\medskipamount}
}}
\begin{document}
\newcommand{\Addresses}{{
  \bigskip
  \footnotesize

}}
\title[On the universal mod $p$ supersingular quotients for $\mathrm {GL}_2(F)$ for a general $F/\Qp$]{On the universal mod $p$ supersingular quotients for $\mathrm {GL}_2(F)$ over $\fpbar$ for a general $F/\Qp$}
\author{{Yotam I. Hendel}}
\begin{center}
\end{center}

\maketitle
\begin{abstract}
Let $F/\Qp$ be a finite extension. 
We explore the universal supersingular mod $p$ representations of $\mathrm{GL}_2(F)$ by computing a basis for their spaces of  invariants under the pro-$p$ Iwahori subgroup. 
This generalizes works of Breuil 
 and Schein 
 (from $\Qp$ and the totally ramified cases to an arbitrary extension $F/\Qp$).
Using these results we then construct, for an unramified $F/\Qp$, a quotient of the universal supersingular module which has as quotients all the supersingular representations of $\mathrm{GL}_2(F)$ with a $\mathrm{GL}_2(\mathcal{O}_F)$-socle that is expected to appear in the mod $p$ local Langlands correspondence.
A construction in the case of an extension of $\Qp$ with inertia degree 2 and suitable ramification index is also presented.
\end{abstract}
\tableofcontents

\pagenumbering{arabic}
\section{Introduction}\label{intro}

Let $p$ be an odd prime, $ F$/$\mathbb Q_p$ a finite extension, $\mathcal {O}_F$ its ring of integers and $\mathbb F_q$=$\mathcal{O}_F/(\pi)$ its 
residue field where $\pi$ is a uniformizer of $\mathcal {O}_F$.
Classifying the smooth irreducible representations of $G=\mathrm{GL}_2(F)$ over $\fpbar$, the algebraic closure of $\mathbb F_q$, is crucial in order to establish a mod $p$ Langlands correspondence. 

The first major progress in this area was done in \cite{BL} by Barthel and Livn\'{e}, and divided the smooth irreducible representations admitting a central character into four possible classes. 
Two key facts in proving this classification are that every irreducible representation of the maximal compact subgroup $K=\mathrm{GL}_2(\mathcal{O}_F)$ arises from an irreducible representation of the finite group $\mathrm{GL}_2(\mathbb F_q)$, and that $\mathrm{End}_G(\indkzg\sigma)\simeq \fpbar[T]$.
Here $\indkzg \sigma$ is the compact induction of an irreducible representation $\sigma$ of $KZ$ (which will be defined in a detailed manner later), $Z$ is the center of $G$ and $\fpbar[T]$ is the polynomial algebra generated by a single element $T$.
Barthel and Livn\'{e} show that in all four cases the representations are quotients of $\indkzg \sigma/(T-\lambda \cdot\mathrm{Id})$ (up to a twist by an unramified character) and describe explicitly three of these classes. 
The remaining fourth class of representations, called supersingular representations, is attained as a quotient of $\indkzg \sigma/(T)$, i.e. for $\lambda=0$, and is still largely a mystery. 
The Barthel-Livn\'{e} classification was later generalized to $\mathrm{GL}_n(F)$ by Herzig \cite{Her},
and it was shown that in order to fully understand the non-supersingular representations of $\mathrm{GL}_n(F)$, one first has to know the supersingular representations for smaller ranks.

If $F\neq \Qp$, even studying supersingular representations of $\mathrm{GL}_2(F)$ is a very hard task \cite{BP},\cite{Schr}. 
The only case in which these representations are fully understood is when $F=\Qp$. 
In this case, the supersingular representations of $G$ were first classified by Breuil  \cite{B} who proved that $\indkzg \sigma/(T)$ are irreducible via computation of the invariant subspace $(\indkzg \sigma/(T))^{I(1)}$, where $I(1)$ is the pro-$p$ Iwahori subgroup, the subgroup of $K$ such that the reduction modulo $p$ of every element is an upper triangular unipotent matrix of $\mathrm{GL}_2(\mathbb F_q)$. 
In \cite{Sc1} and \cite{Sc2}, Schein extended this method to totally ramified field extensions, and proved an irreducibility criterion for supersingular representations.

In this article we use extensions of these methods to explore supersingular representations further, in the case of an arbitrary finite field extension $F/\Qp$.
In particular, we compute the space of $I(1)$-invariants of $\indkzg \sigma / (T)$, as well as construct a universal module of which all the supersingular representations with a prescribed $K$-socle are quotients (recall that the $K$-socle of a $K$-module is the direct sum of its irreducible $K$-submodules).

Considering supersingular representations with a specific $K$-socle is of special interest since it is conjectured that the mod $p$ local Langlands correspondence should give a bijection between certain supersingular representations with a $K$-socle prescribed by generalizations of the weight part of Serre's modularity conjecture (which are now known to be true for $\mathrm{GL}_2(F)$ \cite{GLS}) and certain Galois representations.
This indeed holds for the case $F=\Qp$ where the mod $p$ local Langlands correspondence is completely understood.
\subsection{Overview and main results}
We give an overview of the article and state the main results, we start by establishing the required notations. 
Let $\sigma$ be a smooth irreducible mod $p$ representation of $K$ with $\mathrm{Fr}_p$ the Frobenius automorphism and $\mathrm{Sym}_j^{r_j}\fpbar^2 =\mathrm{Sym}^{r_j} \fpbar^2 \circ \mathrm{Fr}^{p^{j}}_p$. 
Then we have $\sigma = \serwgt$ for some natural numbers $0 \leq r_j \leq p-1$ and $0 \leq \omega \leq q-2$ (see \cite[Proposition 4]{BL}). 
We extend $\sigma$ to $KZ$ by setting the element $\pi \cdot \mathrm{Id} \in Z$ to act trivially.
Note that each $\mathrm{Sym}^{r_j} \fpbar ^2$ can be considered as an $(r_j+1)$-dimensional vector space over $\fpbar$ spanned by $\{x_j^{r_j-k} y_j^k\}_{k=0}^{r_j}$. An elaborate discussion of these representations is given in Section \ref{sec:Not}.

For $v \in V_\sigma$ 
and $g \in G$ we define the function $g \otimes v: G \to V_\sigma$ by  $\sigma(k)(v)$ for $kzg^{-1} \in KZ g^{-1}$ and $0$ elsewhere. 
We then set the {\it compact induction} $\indkzg\sigma$ to be the $\fpbar$-vector space generated by all such functions with the action $g' \cdot (g \otimes v) = g'g\otimes v$.

Now set
$\alpha=
\left(
\begin{array}{ccc} 
1  & 0 \\
0 & \pi \\  \end{array}
\right)$, 
$\beta=\alpha
\left(
\begin{array}{ccc} 
0  & 1 \\
1 & 0 \\  \end{array}
\right)$ and as a set of explicit coset representatives of $\mathcal {O}_F /(\pi^n)$ take $I_{n}=\{\sum\limits_{i=0}^{n-1}\pi^i [\mu_i]:\mu_i \in \mathbb F_q\}$ where $[\mu]$ is the canonical Teichmuller lift of $\mu$ to $\mathcal {O}_F$ and $I_0=\{0\}$.
Recall that the Bruhat-Tits tree of
 $G$ is a $(q+1)$-regular tree with vertices in bijection with the cosets $G/KZ$. 
By the Cartan decomposition we can take explicit representatives:
\[G=\coprod\limits_{n\in\mathbb N}\Bigg(\coprod\limits_{ \mu\in I_n} 
\left(
\begin{array}{ccc} 
\pi ^n & \mu \\
0 & 1 \\  \end{array}
\right)
 KZ\coprod\limits_{\mu\in I_n} 
\beta\left(
\begin{array}{ccc} 
\pi ^n & \mu \\
0 & 1 \\  \end{array}
\right)
\left(
\begin{array}{ccc} 
0 & 1 \\
1 & 0 \\  \end{array}
\right)
 KZ\Bigg),\] and thus every function in $\indkzg\sigma$ is determined by its values on the vertices of the Bruhat-Tits tree of $G$ which are regarded as the cosets
$KZ\left(
\begin{array}{ccc} 
\pi ^n & \mu \\
0 & 1 \\  \end{array}
\right)^{-1}$ and 
$KZ\left(
\begin{array}{ccc} 
1 & 0 \\
\pi \mu & \pi^{n+1} \\  \end{array}
\right)^{-1}$.
We now define explicitly elements of $\indkzg \sigma$ which will be $I(1)$-invariant (mod $(T)$) for specific parameters.
\begin{defn} \label{def:invs2}
For $n\geq1$ and $\mu = \sum\limits_{i=0}^{n-1}\pi^i [\mu_i]$ with $\mu_{i-1} \in \Fq$
 define the following elements of $\indkzg \sigma$:

$\qquad \qquad s^k_n=\sum\limits_{\mu \in I_n}
\left(
\begin{array}{ccc} 
\pi ^n & \mu \\
0 & 1 \\  \end{array}
\right)
\otimes\mu_{n-1}^{k} \xvec$ where $0 \leq k \leq q-1$.

$\qquad \qquad t^s_n=\sum\limits_{\mu \in I_n}{
\left(
\begin{array}{ccc} 
\pi ^n & \mu \\
0 & 1 \\  \end{array}
\right)
\otimes} \xbutsvec \otimes x_s^{r_s-1}y_s$ where $0 \leq s \leq f-1$.
\end{defn} 
Let $I$ denote the Iwahori subgroup of $G$, the following is Theorem \ref{thm:InvBasis}.
\begin{thm} \label{thm:InvBasisIntro} 
Let $e$ and $f$ be the ramification index and inertia degree of $F/\Qp$ respectively, assume $2 < r_l < p-3$ for $0 \leq l \leq f-1$ and define the following sets:

$\qquad\qquad\mathcal S^l_m=\{s_n^{p^l(r_l+1)}\}_{n\geq m}\bigcup\{ \beta s_n^{p^l(r_l+1)}\}_{n\geq m}$, $~~~\mathcal S_m=\bigcup\limits_{l=0}^{f-1} \mathcal S^l_m$,

$\qquad\qquad\mathcal T^l_m=\{t_n^{l}\}_{n\geq m}\bigcup\{ \beta t_n^{l}\}_{n\geq m}$, $~\qquad \qquad \qquad \mathcal T_m=\bigcup\limits_{l=0}^{f-1} \mathcal T^l_m.$

Then an $I$-eigenbasis for the space $(\indkzg\sigma/(T))^{I(1)}$ of $I(1)$-invariants as an $\fpbar$-vector space is given by the images in $\indkzg\sigma/(T)$ of the following sets:
\begin{align*}
&\{\mathrm{Id}\otimes x^r , \alpha \otimes y^r\}  & :e=1, f=1 &\text{ ~(Breuil)}\\
&\mathcal S_2\bigcup \{\mathrm{Id}\otimes x^r , \alpha \otimes y^r\}\bigcup \mathcal T_1  & :e>1, f=1 &\text{ ~(Schein)}\\
&\mathcal S_1\bigcup \{\mathrm{Id}\otimes \xvec , \alpha \otimes \yvec\} & :e=1, f>1\\
&\mathcal S_1\bigcup \{\mathrm{Id}\otimes \xvec , \alpha \otimes \yvec\}\bigcup \mathcal T_1& :e>1, f>1
\end{align*}
\end{thm}
\begin{proof}
Cases $1$ and $2$ were proved by Breuil \cite[Theorem 3.2.4]{B} and Schein \cite[Theorem 2.24]{Sc1} respectively. 
Cases $3$ and $4$ are proved as follows: 
The elements in the sets mentioned above are shown to be $I(1)$-invariant in Proposition \ref{prop:inv1} using computations with the generators of $I(1)$.
The proof is then divided into three parts (Lemmas \ref{lem:InvBasis1}, \ref{lem:InvBasis2} and \ref{lem:InvBasis3}), where it is gradually shown using properties of $T$ (Proposition \ref{prop:inv2}) that an $I(1)$-invariant element in $\indkzg\sigma/(T)$ can be written as a linear combination of the required elements. Finally, using Lemma \ref{lem:char} the elements in $\mathcal S_k$ and $\mathcal T_k$ are indeed distinct since they are $I$-eigenvectors with different eigenvalues (here $I$ is the Iwahori subgroup of $K$).
\end{proof}
\begin{conc}
$\mathrm{End}_G(\indkzg\sigma/(T))=\fpbar$ (for proof see Conclusion \ref{conc:EndScalar}).
\end{conc}
Recall that $K$-socle of a $K$-module is the direct sum of its irreducible $K$-submodules. 
We now wish to use our knowledge of $(\indkzg\sigma/(T))^{I(1)}$ in order to construct a quotient $U_{\bar\rho}$ of $\indkzg\sigma/(T)$ such that every supersingular representation $W$ with a certain $K$-socle that depends on a modular Galois representation $\bar\rho$ factors through. 

Given a suitable irreducible Galois representation $\bar\rho:\mathrm{Gal}(\overline{F}/F) \to \mathrm{GL}_2(\mathbb F_q)$, one can associate to it a multiset $W(\bar\rho)$ of irreducible mod $p$ representations of $K$ via Serre's weight conjecture, 
which for an unramified $F/\Qp$ is parametrized by the subsets of $\{0,\ldots,f-1\}$ 
(for an elaborate discussion of  Serre's weight conjecture and its generalizations see \cite{BLGG}, \cite{BDJ} or \cite{Sc3}). 
The conjectured mod $p$ local Langlands correspondence is expected to associate to $\bar\rho$ a representation 
$\pi(\bar\rho)$ of $G$ whose $K$-socle is exactly $\bigoplus\limits_{\sigma\in W(\bar\rho)}\sigma$ (see \cite{B1}).
As mentioned before, in the established case of $F=\Qp$, it was shown by Breuil that $\mathrm{soc}_K\pi(\bar\rho)=\bigoplus\limits_{\sigma\in W(\bar\rho)} \sigma$ (see \cite{B}). 

In this spirit, we can state the second main theorem of the article, in which we construct explicitly a $G$-module which has as quotients all the supersingular representations of $G$ with a $K$-socle that arises from Serre's weight conjecture for some Galois representation.
For a given generic $\vec r$, this construction holds for an extension $F/\Qp$ where $f=2$ and $e<\min\{\frac{r_j}{2}\}$ (see \ref{sec:GenQuot}).
Note that for the sake of clarity a concrete computation exemplifying the general method for a simple case (quadratic unramified extension of $\Qp$) is presented in Example \ref{ex:QuadExt}.
\begin{thm} \label{thm:UnivQuotIntro}
Let $F/\Qp$ be a finite unramified extension, $\bar\rho: \mathrm{Gal}(\overline {\mathbb Q}_p / F)\to \mathrm{GL}_2(\fpbar)$ be a Galois representation,  and $\sigma_{\varnothing}\in W(\bar\rho)$, then the following holds:
\begin{itemize}
\item[(1)]  There is an explicit construction of a quotient $U_{\bar\rho}$ of $\indkzg \sigma_{\varnothing}$ such that 
$\bigoplus\limits_{\sigma \in W(\bar\rho)} \sigma\subseteq \mathrm{soc}_K(U_{\bar\rho})$ 
and $\mathrm{soc}_K(U_{\bar\rho})$ can only contain 
certain irreducible submodules (See Theorem \ref{thm:UnivQuot} 
and Lemma \ref{lem:Steps1} for the specific stipulations). 
\item[(2)] The quotient $U_{\bar\rho}$ is universal for supersingular representations $W$ of $G$ with $ \mathrm {soc}_K (W)=\bigoplus\limits_{\sigma \in W(\bar\rho)}\sigma$; assume there is a surjective map $\eta:\indkzg\sigma \twoheadrightarrow W$, then there exists a map $\tilde \eta$ such that the following diagram commutes,

$$\begin{tikzpicture}
  \matrix (m) [matrix of math nodes,row sep=3em,column sep=3em,minimum width=2em]
  {
    \indkzg \sigma & &W \\
      & U_{\bar\rho} & \\};
  \path[-stealth]
    (m-1-1) 
            edge [-{>[sep= 2pt]>}] node [above] {$\eta$} (m-1-3) 
	 edge  [-{>[sep= 2pt]>}] node [below] { $\varphi$} (m-2-2)
    (m-2-2) 
            edge [->,dashed] node [below] {$\tilde \eta$} (m-1-3);
\end{tikzpicture}$$
where $\varphi $ is the reduction map taking an element in $\indkzg\sigma$ to its projection in $U_{\bar\rho}$. 
\end{itemize}
\end{thm}
\begin{proof}
The proof utilizes Proposition \ref{prop:invdim}, in which the $K$-submodules of $\indkzg \sigma/(T)$ generated by some of the $I(1)$-invariants that were presented in Theorem \ref{thm:InvBasisIntro} are computed and shown to be irreducible.
Start with $\indkzg\sigma_{\varnothing}$, every $K$-module in $W(\bar\rho)$
can be obtained as a submodule of a quotient (restricted to $K$) of $\indkzg \sigma_{\varnothing}$.
The combinatorial statement that shows this is possible is essentially Lemma \ref{lem:Steps1}, and the procedure is done iteratively by a sequence of quotients by modules of the form $\Phi_\alpha (T_\alpha (\indkzg \sigma_\alpha))$. 
Here, each $\sigma_\alpha$ is an irreducible representation of $K$, $\Phi_\alpha$ is a homomorphism obtained via Frobenius reciprocity from the inclusion of $\sigma_\alpha$ into the previous step and $T_\alpha$ is the generator of $\mathrm{End}_G(\indkzg\sigma_{\alpha})$.
The process depicted above is then done simultaneously to yield a quotient $U$ 
such that $\bigoplus\limits_{\sigma\in W(\bar\rho)} \sigma \subset \mathrm{soc}_K(U)$ and such that $U$ maps surjectively to $W$.
The latter is shown using the fact that $W$ is supersingular, implying that any non-zero map $\eta:\indkzg\tau/(T-\lambda)\to W$ for some $\tau$ must have that  $\lambda = 0$.
In a similar manner, this time using quotients of the form $\Phi_\alpha (\indkzg \sigma_\alpha)$, we produce a quotient $U_{\bar\rho}$ of $U$ such that $f$ factors through $U_{\bar\rho}$ and such that $\mathrm{soc}_K(U_{\bar\rho})$ contains only submodules that participated in the process that produced $U$ (for more details see Theorem \ref{thm:UnivQuot} and Lemma \ref{lem:Steps1}).
\end{proof}

To put this work in context, we mention a few related results.
A work with a similarly explicit approach for an unramified $F/ \Qp$ was carried out by Morra in \cite{Mo2}. 
There, he completely described  the $I$-socle filtration of $\indkzg\sigma/(T)$ where the latter is treated as $K$-module. 
We consider $\indkzg\sigma/(T)$ as a $G$-module, where determining the complete $K$-socle filtration is a harder task, but in turn we gain different information.
Additionally, while we study the spherical Hecke algebra, $\mathrm{End}_G(\indkzg\sigma)$, a different approach is to study modules over the pro-$p$-Iwahori-Hecke algebra $\mathrm{End}_G(\mathrm{ind}_{I(1)}^G 1)$, as done in \cite{Vi1} and \cite{O1}. 
One may pass from a supersingular representation to a representation of the Iwahori-Hecke algebra by taking the $I(1)$-invariants of a supersingular representation, this defines a functor and if $F=\Qp$ this is an equivalence of categories.
The merit of this approach is that the category of representations of the Iwahori-Hecke algebra is much better behaved, but the downside is that if $F\neq\Qp$ we lose information when passing to this category and the $I(1)$-invariants functor ceases to be an equivalence of categories.
\subsection{Acknowledgments}
This work was carried out under the supervision of Dr. Michael Schein in Bar-Ilan University, and the author wishes to express his deepest gratitude and thanks to Michael for being a very knowledgeable, inspiring and patient advisor, and for presenting him this problem.

\subsection{Notations and preliminaries} \label{sec:Not}
Recall that $p$ is an odd prime, $F/\mathbb Q_p$ is a finite extension with $e$ and $f$ the ramification and inertia indices respectively, 
$\mathcal {O}_F$ its ring of integers, 
$\mathbb F_q=\mathcal{O}_F/(\pi)$ its residue field where $\pi$ is a uniformizer of $\mathcal {O}_F$, and $\fpbar$ an algebraic closure of $\mathbb F_q$. 
We also have $G=\mathrm{GL}_2(F)$, where $K=\mathrm{GL}_2(\mathcal{O}_F)\subseteq G$ is a maximal compact subgroup and $Z$ the center of $G$. 
The Iwahori subgroup of $K$, consisting of matrices which are upper triangular modulo $p$, is denoted by $I$ and its pro-$p$-Sylow subgroup by $I(1)$, and we introduce the useful elements 
$ {w}=
\left(
\begin{array}{ccc}
0 & 1 \\
1 & 0 \\  \end{array}
\right)$ and 
$\alpha=
\left(
\begin{array}{ccc} 
1 & 0 \\
0 & \pi \\  \end{array}
\right)$ in $G$.
Note that the element 
 $\beta= \alpha w =
\left(
\begin{array}{ccc} 
0 & 1 \\
\pi & 0 \\  \end{array}
\right)$
normalizes $I(1)$.

The following are some useful results in modular combinatorics which are usually not referenced.
Let $\nu_p(x)$ be the additive $p$-adic valuation of $x\in \mathbb Q$, then we have the following two classical theorems (\cite{SMA},\cite{Fi}): 
\begin{thm}\label{thm:lucas} Set $m,n\in \mathbb N$ and write them in base $p$, i.e. $m=\sum\limits_{s=0}^k m_s p^s$, $n=\sum\limits_{s=0}^k n_s p^s$ and $0\leq n_s,m_s\leq p-1$.
$\\1.$
(Lucas' Theorem)$~$
$\binom{n}{m}\equiv  \prod\limits_{s=0}^k \binom{n_s}{m_s}\mod p$.
$\\2.$ (Legendre's Theorem) Set $s_p(n)=\sum\limits_{s=0}^k n_s$, then
$\nu_p(n!)=\frac{n-s_p(n)}{p-1}$.
\end{thm}
A corollary of Legendre's theorem which will be of use is that for $0<m\leq p^k$, we get $\nu_p(\binom{p^k}{m})=k-\nu_p(m)$. 
\begin{defn}
For positive integers $i\leq r$, write $i=\sum\limits_{j=0}^{k_1}{i_j p^j}$ and $r=\sum\limits_{s=0}^{k_2}{r_j p^j}$ for the $p$-adic expansions of $i$ and $r$. 
We denote by $\vec i$ the vector $(i_0,i_1,\dots,i_{k_1})$ and define $\sum\limits^{\vec{r}}_{\vec{i}=\vec 0}:=\sum\limits_{\vec {i} \in J}$ where $J=\{\vec {i}:0\leq i_j \leq r_j, \forall ~0 \leq j\leq f-1\}$.
We use these notations and the $p$-adic expansion of numbers vastly throughout the paper.
This should cause no confusion.
\end{defn}
Next, set the compact induction of $\sigma$ to be the space of locally constant functions \[\mathrm{ind}_{KZ}^{G} \sigma = \{f:G \to V_{\sigma}| f(hg)=\sigma(h)(f(g)),\forall g\in G, h\in KZ\}\] supported on finitely many cosets of $G/KZ$, with the $G$-action given by $g\cdot f(x)=f(xg)$. 
For $v\in V_{\sigma}$ define
the function $\Id \otimes v \in \indkzg\sigma$ by,
\begin{displaymath}
    (\Id \otimes v)(g) = \left\{
     \begin{array}{lr}
       \sigma(g)(v) & : g\in KZ \\
      0 & : g\notin KZ.
     \end{array}
   \right.
\end{displaymath} 
The space $\mathrm{ind}_{KZ}^{G}\sigma$ is generated as a $G$-module by the set $\{ \mathrm{Id} \otimes v\}_{v \in V_{\sigma}}$.
Another construction which will come into play is $\Ind_I^K \chi$, the usual induction of a representation $\chi$ of $I$ to $K$.

A model for $\indkzg\sigma$ is provided using the Bruhat-Tits tree of $\mathrm{GL}_2(F)$ denoted $\mathcal A$, whose vertices correspond to the cosets $G/KZ$.

Given $n > 0$, consider
$I_{n}=\{[\mu_0]+\pi[\mu_1]+\dots+\pi^{n-1} [\mu_{n-1}]:\mu_i \in \mathbb F_q\}$
where $[\mu]\in \mathcal{O}_F$ is the canonical Teichmuller lift of $\mu \in \Fq$ and set $I_0 = \{ 0\}$.
In order to choose a coordinate system on $\mathcal A$ with which we can do computations easily, for $\mu \in I_n$ we define
$\gn=
\left(
\begin{array}{ccc} 
\pi ^n & \mu \\
0 & 1 \\  \end{array}
\right)$ and 
$\gnn =
\left(
\begin{array}{ccc} 
1 & 0 \\
\pi\mu & \pi^{n+1} \\  \end{array}
\right)$ and note that
 $\alpha=g_{0,0}^1$ and that $\gnn=\beta\gn {w}$.
Given $\mu \in I_m$ where $m>n$ we also set the truncation operators $[\mu]_{n}=\sum\limits_{k=0}^{n-1}\pi^k[\mu_k]$ with $[\mu]_0=0$.

In order to calculate using $\gn$ and $\gnn$, we need to know how to sum Teichmuller representatives. 
This is given in the following lemma.
\begin{lem} \label{lem:WittAdd}
Let $a,b\in \mathbb{F}_q$, then
$[a]+[b]\equiv [a+b]+\pi^e [P_0(a,b)] \mathrm{~mod~} \pi^{e+1}$,
where \[P_0(a,b)= \frac{a^{q^e}+b^{q^e}-(a+b)^{q^e}}{\pi^e}.\]
\end{lem}
\begin{proof}
See \cite[Lemma 2.2]{Sc2}.
\end{proof}
We can now give $G/KZ$ coordinates in term of $\gn$ and $\gnn$ 
using the Cartan decomposition:
\[G=\coprod\limits_{n\in\mathbb N}\Big(\coprod\limits_{ \mu\in I_n} \gn KZ\coprod\limits_{\mu\in I_n} \gnn KZ\Big),\]
implying that an element $\tilde{f}\in \indkzg\sigma$ can be interpreted as a function from $\mathcal A$ to $V_{\sigma}$ which takes values on finitely many vertices.

In the spirit of \cite{B}, we set $S_n^0$ and $S_n^1$ to be the functions supported on cosets of the form $\gn KZ$ and $\gnn KZ$ respectively and also $S_n=S_n^0\coprod S_n^1$, $B_n^*=\coprod\limits_{k=0}^n S_k^*$ and $B_n=\coprod\limits_{k=0}^n S_k$ for $*\in \{0,1\}$. 
The cosets (or vertices) $\gn KZ$ and $\gnn KZ$ are said to be of radius $n$. 

Recalling the representation theory of $K$, if $\sigma:K\to \mathrm{GL}(V_\sigma)$ is an irreducible $\fpbar$-representation, then it arises from a representation of the finite group $\mathrm {GL}_2(\mathbb F_q)$ via inflation. 
A model for $V_{\sigma}$ is given by the space of homogeneous polynomials 
$\mathrm{Sym}^{\vec r}\fpbar^2 = \bigotimes\limits_{j=0}^{f-1}\mathrm{Sym}_j^{ r_j}\fpbar^2$ for some $\vec{r} \in (\mathbb{Z}/p\mathbb{Z})^f$, 
with the basis $\{\bigotimes\limits_{j=0}^{f-1} x_j^{r_j-i_j}y_j^{i_j}\}^{\vec r}_{\vec i=\vec 0}$. 
The action of $\sigma$ is then given by (up to a twist by a power of $\det$):
\[\sigma\left(
\begin{array}{ccc}
a & b \\
 c & d \\  \end{array}
\right  )(\xyvec )=\bigotimes\limits_{j=0}^{f-1} (a^{p^j}x_j+c^{p^j}y_j)^{r_j-i_j}(b^{p^j}x_j+d^{p^j}y_j)^{i_j}.\]
It is easy to see that for such $\sigma$ the space of invariants $\sigma^{I(1)}$ is one dimensional and spanned by the element $\xvec$ and thus if $1\leq r_j\leq p-2$ for $0\leq j \leq f-1$ then $\sigma$ is determined by the $I$ action on $\sigma^{I(1)}$. It is also evident that $\mathrm{Id}\otimes \xvec$ generates $\indkzg\sigma$ as a $G$-module. 

Fix $\overline {\mathbb Q}_p$ an algebraic closure of $\Qp$.
We recall that by Serre's weight conjecture and its generalizations,
 given a continuous irreducible 
Galois representation $\bar\rho: \mathrm{Gal}(\overline {\mathbb Q}_p / F)\to \mathrm{GL}_2(\fpbar)$ 
one can associate to it a multiset $W(\bar\rho)$ of irreducible mod $p$ representations of $\mathrm{GL}_2(\mathbb{F}_q)$. 
For this reason, throughout the article irreducible representations of $\mathrm{GL}_2(\mathbb{F}_q)$ will be called {\it Serre weights}.
This set is of interest to us since it is expected that $\mathrm{soc}_K(\pi(\bar\rho))=\bigoplus\limits_{\sigma\in W(\bar\rho)}\sigma$ where $\pi(\bar\rho)$ is the representation associated to $\bar\rho$ by the conjectural mod $p$ local Langlands correspondence \cite{B1}.
\section{Properties of $I(1)$-invariants}\label{I(1)Comp} 
\subsection{The action of $T$, the generator of $\mathrm{End}_G(\indkzg\sigma)$}
In this section we prove some properties of $I(1)$-invariants that will help us determine $(\indkzg\sigma/T)^{I(1)}$ in Theorem \ref{thm:InvBasis} for a Serre weight $\sigma=\serwgtnew$.
Recall that, 
\[ I(1) =\left\{
\left(
\begin{array}{ccc}
a & b \\
\pi c & d \\  \end{array}
\right) \in I, a\equiv d \equiv 1 \text{ mod } \pi \text{ for}~a,b,c,d\in \mathcal O_F\right\},\]
and notice that the following decomposition holds for every $m\in I(1)$:
\[
\left(
\begin{array}{ccc} 
\pi a+1 & b \\
\pi c & \pi d+1 \\  \end{array}
\right) =\left(
\begin{array}{ccc} 
1 & (\pi d+1)^{-1}b \\
0 & 1 \\  \end{array}
\right)\left(
\begin{array}{ccc} 
1 & 0 \\
\pi c\epsilon^{-1} & 1 \\  \end{array}
\right)\left(
\begin{array}{ccc} 
\epsilon & 0 \\
0 & \pi d+1 \\  \end{array}
\right)\]
where $\epsilon=\pi(a-cb(\pi d+1)^{-1})+1$.
To verify that $s\in (\indkzg\sigma/(T))^{I(1)}$, it is enough to check that $s$ is invariant under the three types of matrices for $a,b,c,d\in \mathcal O_F$,
\[
u^{+}(b) = \left(
\begin{array}{ccc} 
1 &b \\
0 & 1 \\  \end{array}
\right), u^{-}(\pi c)= \left(
\begin{array}{ccc} 
1 & 0 \\
\pi c & 1 \\  \end{array}
\right) \text{and}~u(a,d)=\left(
\begin{array}{ccc} 
\pi a + 1 & 0 \\
0 & \pi d +1 \\  \end{array}
\right).\]
Note that since $Z$ acts trivially the action of the third generator is reduced to the action of
$u(a,0) =\left(
\begin{array}{ccc} 
\pi a +1 & 0 \\
0 & 1 \\  \end{array}
\right)$ where $a \in \mathcal O _F$.

The next proposition gives an explicit description of the operator $T$ which generates $\mathrm {End}_G(\indkzg\sigma)$ and on which further computations rely heavily.

\begin{prop} \label{prop:inv2}  
Let $v=\sum\limits_{\vec i=0}^{\vec r}{c_{\vec i}\xyvec }$. 
For $n\geq 1$ and $ \mu \in I_n$ we have: 

$T(\gn\otimes v) = \sum\limits_{\lambda \in I_1}{g^0_{n+1,\mu+\pi^n\lambda}\otimes\sum\limits^{\vec{r}}_{\vec{i}=0}{c_{\vec i}(-\lambda)^{i}}}\xvec + $
$g^{0}_{n-1,[\mu]_{n-1}}\otimes c_{\vec r}\bigotimes\limits_{j=0}^{f-1}(\mu^{p^j}_{n-1}x_j+y_j)^{r_j}$,

$T(\gnn\otimes v) = \sum\limits_{\lambda \in  I_1}{g^1_{n+1,\mu+\pi^n\lambda}\otimes\sum\limits^{\vec{r}}_{\vec{i}=0}{c_{\vec r - \vec i}(-\lambda)^{i}}}\yvec + g^{1}_{n-1,[\mu]_{n-1}}\otimes c_{\vec 0} \bigotimes\limits_{j=0}^{f-1}(x_j+\mu^{p^j}_{n-1}y_j)^{r_j}$,\\
and for $n=0$,

$T(\mathrm{Id}\otimes v)= \sum\limits_{\lambda \in  I_1}{g^{0}_{1,\lambda}\otimes \sum\limits_{\vec{i}=0}^{\vec{r}}{c_{\vec i}(-\lambda)^{i}}}\xvec+\alpha\otimes c_{\vec r}\yvec$,

$T(\alpha\otimes v)= \sum\limits_{\lambda \in  I_1}{g^{1}_{1,\lambda}\otimes \sum\limits_{\vec{i}=0}^{\vec{r}}{c_{\vec r -\vec i }(-\lambda)^{i}}}\yvec+\mathrm{Id}\otimes c_{\vec 0 }\xvec$.

\end{prop}

\begin{proof}
Using formulas (4) to (8) of \cite{B}, the proof is a straightforward computation as in \cite{Sc2}. 
\end{proof}

\begin{rem}
Note that $T$ is injective since $T(g \otimes v)$ is supported on at least two different neighbors for every vertex $g \in \mathcal A$.
\end{rem}

\subsection{A toolbox for computations regarding $I(1)$-invariants}
\begin{lem} \label{lem:pseriesdesc}
Let $\chi : I \to \fpbar^*$ be a character. 
The principal series representation $\mathrm{Ind}_{I}^{K}{\chi}$ is a $(q+1)$-dimensional vector space and it is generated by $\mathrm{Id} \otimes 1$.
Furthermore, if $\chi$ does not factor through the determinant then it has a unique irreducible quotient and it is of length is $2^f$ as a $K$-module.  
An explicit description for the elements $\tilde{f}\in \mathrm{Ind}_I^K{\chi}$ is given by ($c_{ {w}},c_{\lambda}\in \fpbar$),
\[
\tilde{f}=c_{ {w}} {w} \otimes 1 +
\sum\limits_{\lambda \in I_1}c_{\lambda} 
\left(
\begin{array}{ccc} 
1 & 0 \\
\lambda & 1 \\  \end{array}
\right)
\otimes 1.\]
\end{lem}

\begin{proof}
It can be easily verified that 
$K/I=
\left\{\left(
\begin{array}{ccc} 
0 & 1 \\
1 & 0 \\  \end{array}
\right)I\right\}\cup
\left\{{\left(
\begin{array}{ccc} 
1 & 0 \\
\lambda & 1 \\  \end{array}
\right)I}\right\}_{\lambda \in I_1}$
with $|K/I|=q+1$. 
The description of the socle filtration of the principal series is given in \cite{BS} or in \cite[Theorem 2.4]{BP}.
\end{proof}

\begin{lem} \label{lem:weightConditions}
Let $\sigma=\serwgtnew$ be a Serre weight, $0 < r < q-1$, and let $0\neq v \in \sigma^{I(1)}$. Then
\[\sum\limits_{\lambda \in I_1}\lambda^{q-r-1} 
\left(
\begin{array}{ccc} 
1 & 0 \\
\lambda & 1 \\  \end{array}
\right) 
v  
+(-1)^\omega \left(
\begin{array}{ccc} 
0 & 1 \\
1 & 0 \\  \end{array}
\right) v =0,\]
and we have,

\begin{enumerate}
\item [(i)] $\sum\limits_{\lambda \in I_1} \left(
\begin{array}{ccc} 
1 & 0 \\
\lambda & 1 \\  \end{array}
\right) v =0,$ \qquad \qquad
\item [(ii)] $\sum\limits_{\lambda \in I_1} \left(
\begin{array}{ccc} 
1 & 0 \\
\lambda & 1 \\  \end{array}
\right) v  $
$+(-1)^\omega \left(
\begin{array}{ccc} 
0 & 1 \\
1 & 0 \\  \end{array}
\right) v =0.$
\end{enumerate}
Moreover, if $r=0$, only (i) holds and if $r=q-1$, only (ii) holds.
\end{lem}
\begin{proof}
Every $v\in \sigma^{I(1)}$ is a scalar multiple of $\xvec$, computing:

$\sum\limits_{\lambda \in I_1}\lambda^{q-r-1}\left(
\begin{array}{ccc} 
1 & 0 \\
\lambda & 1 \\  \end{array}
\right)  \xvec = \sum\limits_{\lambda \in I_1}\lambda^{q-r-1}\bigotimes\limits_{j=0}^{f-1}{(x_j+\lambda^{p^j} y_j)^{r_j}}=-\yvec.$

Noting that 
$\yvec=(-1)^\omega \left(
\begin{array}{ccc} 
0 & 1 \\
1 & 0 \\  \end{array}
\right) \xvec$ finishes the proof.

The computations in the case $r=0$ and $r=q-1$ are similar.
\end{proof}

\begin{defn} \label{def:invs}
For the Serre weight $\sigma=\det^w \otimes \mathrm{Sym}^{\vec r}\fpbar^2$ and for $n \geq 0$ we define the following elements of $\indkzg \sigma$:
\begin{enumerate}
\item $t^s_n=\sum\limits_{\mu \in I_n}{
\gn
\otimes} \xbutsvec \otimes x_s^{r_s-1}y_s$ where $0 \leq s \leq f-1$.

\item $s^k_n=\sum\limits_{\mu \in I_n}
\gn
\otimes\mu_{n-1}^{k} \xvec$ where $0 \leq k \leq q-1$, and $s_0^k = \Id \otimes \xvec$.
\end{enumerate}
\end{defn} 
\begin{prop} \label{prop:BasicBasis}
For the Serre weight $\sigma=\serwgtnew$, the set $\{\beta^t s_n^0 : n \in \mathbb {N}_0, t\in \{0,1\}\}$ is a basis for $(\indkzg\sigma)^{I(1)}$. Furthermore, we have
\begin{enumerate}
\item If $0 < r \leq q-1$, the $KZ$-module generated by $s_n^0$ is isomorphic to $\sigma$ and hence irreducible and the one generated by $\beta s_n^0$ is isomorphic to a reducible $(q+1)$-dimensional principal series for all $n\geq0$. 
\item If $r=q-1$, then $s_n^0+\beta s_{n-1}^0$ generates an irreducible one dimensional $KZ$-module 
for all $n\geq1$.
\item If $r=0$, $s_n^0+\beta s_{n-1}^0$ generates an irreducible one dimensional $KZ$-module isomorphic to $\sigma$ for all $n\geq0$. 
\end{enumerate}
\end{prop}

\begin{proof}
The proof is similar to \cite[Proposition 2.5]{Sc1}, we include a proof for the convenience of the reader.

For the first part of the proposition we use \cite[Proposition 14]{BL}. 
We have that $S(G,\mathrm{Sym}^{\vec r}~\fpbar^2)^{I(1)}$ is the module of functions $\varphi:G\to \bigotimes\limits_{j=0}^{f-1}\mathrm{Sym}_j^{r_j}\fpbar^2$ which are locally constant, compactly supported modulo $KZ$ 
on the left and for which $\varphi(kgi)=\sigma(k)\varphi(g)$ for all $k\in KZ$ with $g\in G$ and $i\in I(1)$. 
The proposition states that $S(G,\mathrm{Sym}^{\vec r}~\fpbar^2)^{I(1)}$ has a basis of $\{\psi_n\}_{n\in \mathbb Z}$ where each $\psi_n$ is supported on $KZ\alpha^{-n}I(1)$ and satisfies $\psi_n(\alpha^{-n})=\xvec $ if $n\leq 0$ and $\psi_n(\alpha^{-n})=\yvec $ otherwise. 
Also, we have $\coprod\limits_{\mu\in I_n} \gn KZ=I(1)Z\alpha^{-n}KZ$ and $\coprod\limits_{\mu \in I_n} \gnn KZ=I(1)Z\beta\alpha^{-n}KZ$ (cf. \cite{BL}).
Since $g\otimes v\in \indkzg\sigma$ represents a function supported on $KZg^{-1}$, we see that $\psi_n=A_n^0(\sigma)$ and $\psi_{-n-1}=A_n^1(\sigma)$, proving the first part of the proposition. 

Note that $A_n^1(\sigma)$ and $A_n^0(\sigma)$ are eigenvectors for the standard $I$-action \cite[Proposition 15]{BL}. 
In order to compute the $KZ$-modules they generate it is enough to compute their linear combinations over $KZ/I$, as given in Lemma \ref{lem:weightConditions}:
$$KZ/I=
\left\{Z \left(
\begin{array}{ccc} 
0 & 1 \\
1 & 0 \\  \end{array}
\right)I\right\}\cup
\left\{{Z \left(
\begin{array}{ccc} 
1 & 0 \\
\lambda & 1 \\  \end{array}
\right)I}\right\}_{\lambda \in I_1}.$$
The $KZ$-module generated by $A_0^0(\sigma)=\mathrm{Id}\otimes \xvec$ is irreducible and isomorphic to $\sigma$. 
Computing the $KZ$ linear combinations of ${A_1^0(\sigma)}$, we see that each is supported on a different coset of the same radius, and thus they are linearly independent and the generated module is $q+1$ dimensional.

Since $A_0^1(\sigma)$ is an $I$-eigenvector, we have a $\varphi\in \mathrm {Hom}_I(\chi_{A_0^1},\indkzg\sigma_{|I})$ for the suitable character $\chi_{A_0^1}$ of $I$ such that $\varphi(1)=(-1)^w\alpha\otimes \yvec$. 
Using Frobenius reciprocity \cite[Subsection 2.1]{BL}, we get a map of $K$-modules $B(\varphi):\mathrm{Ind}_I^K\chi_{A_0^1}\to\indkzg\sigma$ and $B(\varphi)(\mathrm{Id}\otimes1)=\alpha\otimes\yvec$.

Since $\mathrm{Id}\otimes 1$ generates $\mathrm{Ind}_I^K\chi_{A_0^1}$, the image of $\mathrm{Ind}_I^K\chi_{A_0^1}$ is generated by 
 $\alpha\otimes\yvec$. 
As $\alpha\otimes\yvec$ generates a $q+1$ dimensional module, $B(\varphi)$ is an isomorphism onto its image by dimension considerations, and thus $A_{0}^1(\sigma)$ generates a $q+1$-dimensional principal series.
Recalling that $T$ is an injective map of $G$-modules and that $T(\An)=A_{n+1}^0(\sigma)$ and $T(\Ann)=A_{n+1}^1(\sigma)$ the proof is complete.

Let $r=q-1$, then $A_{1}^0(\sigma)$ and $A_{0}^1(\sigma)$ are $I(1)$ invariants and so are their sum.
Since $\lambda^{r}=\lambda^{q-1}=1$, we see that $A_{1}^0(\sigma)+A_{0}^1(\sigma)$ is an eigenvector for the action of the cosets of $K/I$:
\begin{flalign*}
\left(
\begin{array}{ccc} 
1 & 0 \\
\lambda & 1 \\  \end{array}
\right) (A_{1}^0(\sigma)+A_{0}^1(\sigma)) 
=& \sum\limits_{\nu \in I_1}{{g}^0_{1,\nu}\otimes  (1-\lambda\nu)^{r}\xvec} \\
&+ \alpha \otimes \lambda^{r}\yvec 
\left(
\begin{array}{ccc} 
0 & 1 \\
1 & 0 \\  \end{array}
\right) (A_{1}^0(\sigma)+A_{0}^1(\sigma)) \\
=& \sum\limits_{\nu \in I_1}{{g}^0_{1,\nu}\otimes  (-1)^{w}(-\nu)^{r}\xvec} \\
&+g^0_{1,0}\otimes (-1)^w \xvec+ (-1)^{w}\alpha \otimes \yvec.
\end{flalign*}
Using the injectivity of $T$ as before, one gets that $\An+A^1_{n-1}$ generates an irreducible one dimensional module.

Let $r=0$, and observe that 
$T(A_{0}^0(\sigma))= \sum\limits_{\lambda \in I_1}{g^{0}_{1,\lambda}\otimes} 1+\alpha\otimes 1=A_{1}^0(\sigma)+A_{0}^1(\sigma)$.
Now, $A_0^0(\sigma)=\mathrm{Id}\otimes 1$ generates a $1$-dimensional $KZ$-subomodule isomorphic to $\sigma$, and so does its image under $T$.
Since $T^2(A_0^0(\sigma))=A_2^0(\sigma)+A^1_1(\sigma)+A_0^0(\sigma)$, for $m\geq 0$ we have that
 $T^m(A_{0}^0(\sigma))=\sum\limits_{s=0}^{m}{\sum\limits_{\mu \in I_s}{g_{s,\mu}^{1-\delta_{s,m}}\otimes 1}}$ where $\delta_{s,m}$ is the parity Kronecker delta function taking $1$ when $s$ and $m$ are of the same parity and $0$ otherwise.
Hence, $T^n(A_{0}^0(\sigma))-T^{n-2}(A_{0}^0(\sigma))=A_{n}^0(\sigma)+A_{n-1}^1(\sigma)$ for all $n\geq2$.
\end{proof}

The next proposition generalizes \cite[Corollary 2.6]{Sc1}, and we will use it extensively.
\begin{prop}\label{prop:wgtgen1}
Let $W$ be a $G$-module, $x\in W$ and assume that $x$ is an $I(1)$-invariant and that the $KZ$-submodule of $W$ generated by $x$ is irreducible and isomorphic to a Serre weight $\sigma=\serwgtnew$.
Define $x_n^0=\sum\limits_{\mu\in I_n}g_{n,\mu}^0 x$, $x_n^1=\beta x_n^0 =\sum\limits_{\mu\in I_n}g_{n,\mu}^1  {w} x$.

Then for $n\geq0$ the elements $x_n^s$ for $s=0,1$ are $I(1)$-invariant. 
If $r > 0$, then for $n\geq 0$ if $x_n^0\neq 0$ it generates a $KZ$-submodule of $W$ which is irreducible and isomorphic to $\sigma$.
If $r=0$ then for $n\geq0$ the $KZ$-submodule of $W$ generated by $x_n^0+(-1)^w x_{n-1}^1$ is either trivial or isomorphic to $\sigma$.
\end{prop}

\begin{proof}
The proof is analogous to the proof of  \cite[Corollary 2.6]{Sc1}.
%
\end{proof}
\begin{lem} \label{lem:poly}
Let $n \geq 1$.
For every set-theoretic map $f:I_n \to \fpbar$ there exists a unique polynomial $P \in \fpbar[z_0,..,z_{n-1}]$ in which each variable appears with degree at most $q-1$ 
and such that $f(\mu)=P(\mu_0,\mu_1,\dots,\mu_{n-1})$ for all $\mu\in I_n$.
\end{lem}
\begin{proof}
For $n=1$ the argument of \cite[Lemma 3.1.6]{B} can be used with a slight change, the result then follows from the proof of \cite[Lemma 2.1]{Sc2}.
\end{proof}

\section{Computation of $I(1)$-invariants} \label{sec:3}
\subsection{$I(1)$-invariants in $\indkzg\sigma/(T)$}
In this section we will compute explicitly the $I(1)$-invariant space of $\indkzg\sigma /(T)$ and the $KZ$-modules generated by the elements of this space.
Note that henceforth $r= \sum\limits_{j=0}^{f-1}r_j p^j$ denotes the parameter arising from the representation $\sigma=\serwgtnew$ of $K$ and we assume that $2<r_j<p-3$ for $0 \leq j \leq f-1$ unless stated otherwise.
We start by noting that for $a,b,c \in \mathcal{O}_F$ and $[\mu_0]=\mu \in I_1$ the following identities hold: 
\[\left(
\begin{array}{ccc} 
1 & b \\
0 & 1 \\  \end{array}
\right)
\left(
\begin{array}{ccc} 
\pi  & \mu \\
0 & 1 \\  \end{array}
\right) 
=\left(
\begin{array}{ccc} 
\pi & [b_0+\mu_0] \\
0 & 1 \\  \end{array}
\right)
\underbrace{\left(
\begin{array}{ccc} 
1 & B(\mu,b) \\
0 & 1 \\  \end{array}\right)}_\emph{Acts trivially on $\xvec$.}, \]
where $B(\mu,b)=\pi ^{e-1} [P_0(\mu,b_0)]+[b_1]+\pi [b_2] +\pi^2 [b_3]+\ldots$
and $P_0(\mu,b_0)$ is a polynomial arising from addition in the ring of Witt vectors as in Lemma \ref{lem:WittAdd}.
\[\left(
\begin{array}{ccc} 
1 & 0 \\
\pi c & 1 \\  \end{array}
\right)
\left(
\begin{array}{ccc} 
\pi  & \mu \\
0 & 1 \\  \end{array}
\right) 
=\left(
\begin{array}{ccc} 
\pi & \mu \\
0 & 1 \\  \end{array}
\right)
\underbrace{\left(
\begin{array}{ccc} 
1-\pi c\mu  & -\mu^2 c \\
\pi^2 c & 1+\pi c \mu \\  \end{array}\right)}_\emph{Acts trivially on $\xvec$.}, \]
\[\left(
\begin{array}{ccc} 
\pi a+1 & 0 \\
0 & 1 \\  \end{array}
\right)
\left(
\begin{array}{ccc} 
\pi  & \mu \\
0 & 1 \\  \end{array}
\right) 
=\left(
\begin{array}{ccc} 
\pi & \mu \\
0 & 1 \\  \end{array}
\right)
\underbrace{\left(
\begin{array}{ccc} 
\pi a +1 & a \mu \\
0 & 1 \\  \end{array}\right)}_\emph{Acts trivially on $\xvec$.} .\]
\begin{defn}
Let $M \leq \indkzg \sigma/(T)$ be an $H$-submodule where $I(1) \leq H \leq G$. 
We say that an element $s \in \indkzg \sigma/(T)$ is a non-trivial $I(1)$-invariant$\mod M$ if its image in $(\indkzg \sigma/(T))_{|H}/M$ is a non-trivial $I(1)$-invariant.
\end{defn}
\begin{prop} \label{prop:inv1} 
Let $s_n^k$ and $t_n^s$ be the elements of $\indkzg\sigma/(T)$ as defined in Definition \ref{def:invs}, then we have the following: 
\begin{enumerate}
\item If $n\geq1$ and $k \neq r$ we have $s^k_n\in \mathrm{Im}(T)$ for $0\le k_j \le r_j$ and $ 0 \leq j \leq f-1$.

\item If $n=1$ the element $s_n^r$  is a non-trivial $I(1)$-invariant. 
If  $n\geq2$, then $s_n^r\in \mathrm{Im}(T)$.

\item 
If $f>1$, for $n \geq 1$ and $0\leq l \leq f-1$ 
then $s_n^{p^l (r_l+t)}$ is a non-trivial $I(1)$-invariant$\mod <\{{s_n^{p^l (r_l+s)}\}_{0 \leq s \leq t-1}{}}>_{G}$ where $0\le t \le p-r_l-1$. 
If $f=1$ the claim is true for $n \geq 2$.

\item 
Set $n\geq1$. 
If $e>1$, for $0\le k \le f-1$ and  $r>2p^s$ the element $t_n^s$ is a non-trivial $I(1)$-invariant. 
If $e=1$, then $t_n^s$ is a non-trivial $I(1)$ invariant $\mod<\{s_n^{p^{(f+s-1)}m}\}_{1\leq m \leq p-1} >_{G}$.
\end{enumerate}
\end{prop}

\begin{proof}~
\begin{enumerate}
\item For $n=1$, by Proposition \ref{prop:inv2}, taking $v=\bigotimes \limits _{j=0}^{f-1} x_j^{r_j-k_j}y_j^{k_j}$ with 
 $k\neq r$ and $c_{\vec k}=1$, we have
\[T(\mathrm{Id}\otimes v)=
\sum\limits_{\mu \in I_1}{g^{0}_{1,\mu}\otimes(-\mu_0)^{k}}\xvec = (-1)^k s_1^k.\]
For $n \geq 2$, since $T$ is $G$-equivariant we compute to get that,
\[T((-1)^k \sum\limits_{\mu \in I_{n-1}}g^0_{n-1,\mu}  \otimes v)
=(-1)^k\sum\limits_{\mu \in I_{n-1}}g^0_{n-1,\mu} T(\Id \otimes v)= \sum\limits_{\mu \in I_{n-1}}g^0_{n-1,\mu} s_1^k=s_n^k.\]

\item 
Notice that by Proposition \ref{prop:inv2} we have
\[
T(\mathrm{Id} \otimes \yvec ) = (-1)^r s_1^r+ \alpha \otimes \yvec.
\]
Since $\beta$ normalizes $I(1)$, the element $\alpha \otimes \yvec=\beta (\mathrm{Id}\otimes \xvec)$ is an $I(1)$-invariant and thus so is $s_1^r$. 
Furthermore, $s_1^r \notin \mathrm{Im}(T)$ since $\mathrm{Id}\otimes \xvec \notin \mathrm{Im}(T)$.

If $n\geq 2$ we see that $s_n^r \in \mathrm{Im}(T)$ since
\[T(\sum\limits _{\mu \in I_{n-1}}g^0_{n-1,\mu}\otimes \yvec) =
(-1)^r s_n^r+\sum\limits_{\mu\in I_{n-2}}g^0_{n-2,[\mu]_{n-2}}\otimes \sum\limits_{\lambda \in I_1} \bigotimes\limits_{j=0}^{f-1}(\lambda^{p^j}x_j+y_j)^{r_j},\]
where the right-most expression sums up to $0$ because $\sum\limits_{\lambda \in I_1} \lambda^k=0$ for $k \neq q-1$.
\item
Considering elements of the form $s_n^{p^l(r_l+ t)}$,
 we see in view of Proposition \ref{prop:inv2} that the coefficient of $\xvec$ is polynomial in $\mu_{n-1}$ of degree that cannot appear in $\mathrm{Im}(T)$.  By Lemma \ref{lem:poly}, these polynomials have a unique presentation so $s_n^{p^l(r_l+ t)}$ are not trivial.
Take $n=1$ and $0 \leq t \leq p-r_l -1$, we start by computing the action of the three generator types as the remark in the beginning of Section $2$ suggests. 
The elements $u^{-}(\pi c)$ and $u(a,0)$ act trivially on $s_1^k$ and in particular
$u^-(\pi c)u(a,0)(s_1^{p^l(r_l+t)}) -s_1^{p^l(r_l+t)}=0$ for all $a,c \in \mathcal O_F$.
Applying $u^{+}(b)$,
\begin{flalign*}
\left(
\begin{array}{ccc} 
1 & b \\
0 & 1 \\  \end{array}
\right)
s_1^{p^l(r_l+t)} -s_1^{p^l(r_l+t)}&=\sum\limits_{\mu \in I_1}{g}^0_{1,\mu}\otimes\Big( (\mu-b_0)^{p^l(r_l+t)} -\mu^{p^l(r_l+t)}\Big) \xvec \\
&=\sum\limits_{s=0}^{r_l-1}(-b_0)^{p^l(r_l+t-s)}\binom{r_l+t}{s_l}s_1^{p^l s} \\
&+\sum\limits_{m=0}^{t-1}(-b_0)^{p^l (t-m)} \binom{r_l+t}{r_l+m} s_1^{p^l(r_l+m)}.
\end{flalign*}
Using (1.), we see that $s^{p^l s}_1\in \mathrm{Im}(T)$ for $0 \leq s \leq r$, and recalling the assumptions we see that $s_1^{p^l(r_l+t)}$ is invariant as claimed.
If $f=1$, then $s_1^{p^l r_l}=s_1^r$, which is not trivial, and in that case the elements are invariant for $n\geq 2$.
Using Proposition \ref{prop:wgtgen1} with $W=(\indkzg\sigma/(T))/<\{{s_n^{p^l(r_l+s)}\}_{0 \leq s \leq t-1}{}}>_{G}$ and $x=s_1^{p^l (r_l+t)}$ finishes the proof.
\item 
Recalling $e$ is the ramification index of $F/\mathbb Q_p$ (and $\delta_{e,h}$ the Kronecker delta function),
\begin{flalign*}
\left(
\begin{array}{ccc} 
1 & b \\
0 & 1 \\  \end{array}
\right)t_1^s -t_1^s
&=\sum\limits_{\mu \in I_1}{g}^0_{1,[\mu+b_0]}\otimes \xbutsvec \otimes x_s^{r_s-1}(B(\mu,b)^{p^s}x_s+y_s)-t_1^s\\
&=\sum\limits_{\mu \in I_1}{g}^0_{1,\mu}\otimes \xbutsvec \otimes x_s^{r_s-1}((\delta_{e,1} P_0(b,[\mu-b_0])+b_1)^{p^s}x_s+y_s)-t_1^s
\\
&=\delta_{e,1}\sum\limits_{\mu \in I_1}{g}^0_{1,\mu}P_0(b_0,\mu-b_0)^{p^s} \otimes \xvec+\underbrace {b_1^{p^s}s_1^0}_{\emph{$\in \mathrm{Im}(T)$}}.\\
\left(
\begin{array}{ccc} 
1 & 0 \\
\pi c & 1 \\  \end{array}
\right)
t_1^s -t_1^s
&=\sum\limits_{\mu \in I_1}{g}^0_{1,[\mu]}(-c\mu ^2 )^{p^s} \otimes \xvec = (-c)^{p^s}s_1^{2p^s}.\\
\left(
\begin{array}{ccc} 
\pi a +1 & 0 \\
0 & 1 \\  \end{array}
\right)
t_1^s- t_1^s
&=\sum\limits_{\mu \in I_1}{g}^0_{1,[\mu]}(a \mu)^{p^s} \otimes \xvec = a^{p^s}s_1^{p^s}.
\end{flalign*}

We see that $t_1^s$ is clearly an $I(1)$-invariant if $r>2p^s$ and $e>1$. 
If $e=1$, the degrees of $\mu$ appearing in $P_0(b_0,\mu-b_0)^{p^s}$ are $\{p^{(s-1)}m\}_{1\leq m \leq p-1}$.
We get, \[P_0(b_0,\mu-b_0)^{p^s}
=\sum\limits_{k=1}^{p-1}\mu^{kp^{f+s-1}} (-b_0)^{p^{f+s-1}(p-k)} \frac{\binom{p}{k}}{p}
=\sum\limits_{k=1}^{p-1} c_k(b_0,s) \mu^{kp^{f+s-1}},\] 
for some non-zero constants $c_s(b_0,s)$. 
Thus, $t_1^s$ is an invariant$\mod<\{s_1^{p^{(s-1)}m}\}_{1\leq m \leq p-1}>_{G}$.
Furthermore, $t_1^s$ is not trivial as it clearly is not in the image of $T$. 
The result for arbitrary $n\geq 2$ follows using Proposition \ref{prop:wgtgen1}.
\end{enumerate}
\end{proof}
\begin{cor} \label{cor:lower degrees are trivial}
Let $n \geq 1$, and let $Q(z_0,\ldots , z_{n-1}) \in \fpbar[z_0,\ldots,z_{n-1}]$ be a polynomial in which each variable 
appears with degree at most $q-1$.
Set $d_{z_{n-1}} :=\deg_{z_{n-1}}(Q)$, and assume $d_{z_{n-1}} \neq r$ and that each digit $(d_{z_{n-1}})_j$ in the $p$-adic expansion of $d_{z_{n-1}}$ satisfies $0 \leq (d_{z_{n-1}})_j \leq r_j$.
Then the element 
$\sum\limits_{\mu \in I_n}g^{0}_{n,\mu}\otimes Q(\mu)\xvec $ lies in $\mathrm{Im}(T)$.
\end{cor}
\begin{proof}
By Proposition \ref{prop:inv1}(1) we have $s^k_1 \in \mathrm{Im}(T)$ if $0 \leq k_j \leq r_j$ and $k \neq r$, which covers the case $n=1$.
If $n \geq 2$, let $Q'$ be the polynomial obtained from $Q$ by substituting $z_{n-1}=1$.
The claim now follows since for every $k$ as above,
\[
\sum\limits_{\mu \in I_{n}}g^{0}_{n,\mu} \otimes Q'([\mu]_{n-1})\mu_{n-1}^k  \xvec
=\sum\limits_{\mu \in I_{n-1}}g^{0}_{n-1,\mu} Q'(\mu)s_1^k \in \mathrm{Im}(T).
\]
\end{proof}
The following example shows that the elements in Proposition \ref{prop:inv1} are the correct generalization of the elements $Q_n^0(\sigma)$ and $Q_n^1(\sigma)$ defined in \cite{Sc1}.
Also, the computation itself shows that invariants arising from different digits, i.e. elements $t_n^s$  and $s_n^{p^s(r_s+1)}$ for different $0 \leq s \leq f-1$, are somewhat independent.
\begin{exmpl}
Take $n \geq 1$, if $f>1$ and $r_j>0$ for $0 \leq j \leq f-1$, then the elements $s_n^{r+p^lt}$ are not $I(1)$-invariant$\mod<\{s_n^{r+mp^l}\}_{1\leq m \leq t-1} >_{G}$ for $0 \leq t\leq p-r-1$ and $0 \leq l \leq f-1$.
\end{exmpl}
\begin{proof}
If $f=1$, then $s_n^{r+p^lt}=s_n^{r_0+p^lt}$ and this is an invariant for $n\geq 2$ by Proposition \ref{prop:inv1} or by \cite{Sc1}.
Otherwise, calculate (recall that we write $m_j$ for the $p$-adic expansion of $m$),
\[
\left(
\begin{array}{ccc} 
1 & b \\
0 & 1 \\  \end{array}
\right)
s_1^{r+p^l t} -s_1^{r+p^l t}
=\sum\limits_{m=0}^{r+p^l t-1}(-b_0)^{r+p^l t-m}\bigg (\prod\limits_{l\neq j=0}^{f-1}\binom{r_j}{m_j}\bigg) \binom{r_l+t}{m_l} s_1^m.\]
But elements of the form $s_1^{r -p^k+t p^l}$ for some $k \neq l$ exist in this sum and cannot be obtained via $T$.
Since $u^{+}(b) \sum\limits_{\mu \in I _1}g^0_{1,\mu} s_n^{k}=\sum\limits_{\mu \in I _1}g^0_{1,\mu+b_0}u^+({B(\mu,b)}) s_n^k $ for the $B(\mu,b)$ presented in Proposition \ref{prop:inv1}, we see that $s_n^{r+p^l t}$ cannot be invariant for $n \geq 2$.
\end{proof}
\subsection{The Serre weights generated by $I(1)$-invariants in $\indkzg\sigma /(T)$}
We turn to develop tools in order to compute the $KZ$-modules generated by different $I(1)$-invariants.
\begin{lem} \label{lem:FinInd}
Assume $M$ is a $K$-module and $s\in M$ is an $I(1)$-invariant on which $I$ acts via a character $\chi$, where $\chi$ does not factor through the determinant and $\chi(\pi \Id)=1$. Define 
\[ \ell_{s} =\Big\{x
\in \mathrm{Ind}_I^K\chi: 
c_{ {w}} {w} s +
\sum\limits_{\lambda \in I_1}c_{\lambda}\left(
\begin{array}{ccc} 
1 & 0 \\
{\lambda} & 1 \\  \end{array}
\right)s =0  
\Big \},\]
where we use the description of $x$ given in Lemma \ref{lem:pseriesdesc}, then $\mathrm{Ind}_I^K \chi/\ell_s \simeq<s>_K$. 
Furthermore, if $\dim<s>_{KZ}\leq \dim \sigma_{\chi}$, then $<s>_{KZ}\simeq \sigma_{\chi}$ where $ \sigma_{\chi}$ is the unique Serre weight where the $I$-action on $\sigma_\chi^{I(1)}$ is via $\chi$.
\end{lem}
\begin{proof}
Let $s$ be an $I(1)$-invariant and assume that $I$ acts on $s$ via an $\fpbar$-valued character $\chi$. 
This defines a map $\varphi\in \mathrm {Hom}_I(\chi,M_{|I})$, with $\varphi(1)=s$.
Consider the map $B(\varphi):\mathrm{Ind}_I^K\chi\to M$ obtained by Frobenius reciprocity and note that $B(\varphi)(\mathrm{Id}\otimes 1)=s$.
Since $\mathrm{Ind}_I^K\chi$ is spanned by $\mathrm{Id}\otimes 1$ as a $K$-module, 
$<{B(\varphi)(\mathrm{Id}\otimes 1)}>_K=\mathrm{Im} B(\varphi)=<{s}>_K\subseteq M$, 
and since $\mathrm{Im} B(\varphi) \simeq \mathrm{Ind}_I^K \chi/ \ker B(\varphi) $, in order to calculate $<{s}>_{KZ}$, one can compute $\ker B(\varphi)$.

By Lemma \ref{lem:pseriesdesc} every $f \in \Ind_I^K \chi$ is of the form
$f=c_{ {w}} {w}\otimes 1  +
\sum\limits_{\lambda \in I_1}c_{\lambda}\left(
\begin{array}{ccc} 
1 & 0 \\
{\lambda} & 1 \\  \end{array}
\right)\otimes 1$, and since $B(\varphi)(\mathrm{Id}\otimes 1)=s$, by the following $\ker(B(\varphi))=\ell_s$,
\[
B(\varphi)(f)=c_{ {w}} {w}s  +
\sum\limits_{\lambda \in I_1}c_{\lambda}\left(
\begin{array}{ccc} 
1 & 0 \\
{\lambda} & 1 \\  \end{array}
\right)s.
\]

For the second part, since $\mathrm{Ind}_I^K\chi$ has a unique irreducible quotient by Lemma \ref{lem:pseriesdesc}, so does
$<s>_K \simeq \mathrm{Ind}_I^K \chi/\ell_s$, implying that $\dim \sigma_\chi \leq \dim <s>_{KZ}$.
\end{proof}
\begin{rem}
Since $I(1)$ is pro-$p$ every character $\chi: I(1) \to \fpbar^\times$ is trivial.
In particular, this means that every $I$-eigenvector $\tilde{f} \in \indkzg \sigma$ is an $I(1)$-invariant (see \cite[Lemma 3]{BL}).
\end{rem}
Equipped with Lemma \ref{lem:FinInd}, we turn to compute the dimensions of the $KZ$-submodules generated by the invariants found earlier in this section. 
We begin with a computation for the representatives of $K/I$:
\begin{enumerate}
\item $
\left(
\begin{array}{ccc} 
0 & 1 \\
1 & 0 \\  \end{array}
\right){g}^0_{1,\mu}= {g}^0_{1,\mu^{-1}}
\left(
\begin{array}{ccc} 
-\mu^{-1} & 0 \\
0 & \mu \\  \end{array}
\right)\left(
\begin{array}{ccc} 
1 & 0 \\
\mu^{-1} \pi & 1 \\  \end{array}
\right)$ for
$\mu \neq 0$.

\item $
\left(
\begin{array}{ccc} 
0 & 1 \\
1 & 0 \\  \end{array}
\right){g}^0_{1,\mu}=\alpha
\left(
\begin{array}{ccc} 
0 & 1 \\
1 & 0 \\  \end{array}
\right)$ for $\mu=0$.
\item Let $\nu = \mu(1+\lambda \mu)^{-1}$ and assume $\mu\neq-\lambda ^{-1}$,\\
$
\left(
\begin{array}{ccc} 
1 & 0 \\
\lambda & 1 \\  \end{array}
\right){g}^0_{1,\mu}=
{g}^0_{1,\nu} \left(
\begin{array}{ccc} 
(1+\lambda\mu)^{-1} & 0 \\
\lambda \pi & 1+\lambda \mu \\  \end{array}
\right)$
$={g}^0_{1,\nu}
\left(
\begin{array}{ccc} 
1-\lambda \nu & 0 \\
\lambda \pi & (1-\lambda \nu)^{-1} \\  \end{array}
\right)$.
\item $
\left(
\begin{array}{ccc} 
1 & 0 \\
\lambda & 1 \\  \end{array}
\right){g}^0_{1,\mu}=\alpha
\left(
\begin{array}{ccc} 
1 & \lambda^{-1}\pi \\
0 & 1 \\  \end{array}
\right)
\left(
\begin{array}{ccc} 
-\lambda^{-1} & 0 \\
0 & \lambda \\  \end{array}
\right)
\left(
\begin{array}{ccc} 
0 & 1 \\
1 & 0 \\  \end{array}
\right)
$ for $\mu=-\lambda ^{-1}$.
\end{enumerate}
\begin{lem}\label{lem:char}
Let $\left(
\begin{array}{ccc} 
a & b \\
c & d \\  \end{array}
\right)=i \in I$, 
then its action on the invariants as in Proposition \ref{prop:inv1} is the following:

$i)$
$i\cdot s_n^k=a^{r-2k}(ad)^{k+\omega} s_n^k$. Note that $r-2k=\sum\limits_{j=0}^{f-1}(r_j-2k_j)p^j$.

$ii)$
$i\cdot t_n^s=a^{r-2p^s}(ad)^{p^s+\omega}t_n^s$.  Note that $r-2p^s=\sum\limits_{j=0}^{f-1}(r_j-2 \delta_{j,s})p^j$.

\end{lem}

\begin{proof}
The $I$-action factors through $I/I(1)=
\bigg\{\left(
\begin{array}{ccc} 
a & 0 \\
0 & d \\  \end{array}
\right)\bigg\}_{a,d \in \Fq^\times}
$:

$i)$
$i\cdot s_n^k= \sum\limits_{\mu \in I_n}g_{n,a d^{-1}\mu}^0\otimes \mu_{n-1}^k a^{r}(ad)^{\omega}\xvec
=a^{r-2k}(ad)^{k+\omega} s_n^k$.

$ii)$
$i\cdot t_n^s= \sum\limits_{\mu \in I_n}g_{n,a d^{-1}\mu}^0\otimes  a^{r-p^s}d^{p^s}(ad)^{\omega} \xbutsvec\otimes x_s^{r_s-1}y_s
=a^{r-2p^s}(ad)^{p^s+\omega}t_n^s$.
\end{proof}

\begin{prop} \label{prop:invdim}
\renewcommand{\theenumi}{\roman{enumi}}
The images of the following elements, when they are $I(1)$-invariant, generate in the proper quotients of $\indkzg\sigma/(T)$ irreducible $KZ$-submodules of the following dimensions (where the exact conditions and quotients are as given in Proposition \ref{prop:inv1}):
\begin{enumerate}
\item The element $\mathrm{Id}\otimes \xvec$ generates an irreducible $KZ$-submodule isomorphic to $\sigma$.
\item If $f>1$, then $s_n^{p^k(r_k+1)}$ generates an irreducible $KZ$-submodule isomorphic to \\
$\det^{\omega+p^k(r_k+1)} \otimes\mathrm{Sym}^{\vec r-2p^k(r_k+1)}\fpbar^2$
for all $n \geq 1$. 
If $f=1$ it generates  an irreducible $KZ$-submodule isomorphic to 
$\det^{\omega+r+1} \otimes\mathrm{Sym}^{p-r-3}\fpbar^2$ for all $n\geq2$.
\item The element $t_n^s$ generates an irreducible $KZ$-submodule isomorphic to 
$\det^{\omega+p^s} \otimes\mathrm{Sym}^{\vec r-2p^s}\fpbar^2$
for all $n \geq 1$.
\end{enumerate}
\end{prop}
\begin{proof}~
\renewcommand{\theenumi}{\roman{enumi}}
\begin{enumerate}
\item  
We have the quotient map $\indkzg \sigma \to \indkzg \sigma/ (T)$ sending $\mathrm{Id}\otimes \xvec$ to its image.
By Frobenius reciprocity we have a non-zero map $\psi:\sigma \to \indkzg \sigma/ (T)_{|KZ}$ sending $\xvec$ to $\mathrm{Id}\otimes \xvec \in \indkzg \sigma/ (T)_{|KZ}$. Since $\sigma$ is irreducible $\psi$ is an isomorphism onto its image.
\item
Assume $f>1$ and set $k'=p^k(r_k+1)$ for convenience. 
In a fashion similar to previous calculations, one computes ($\lambda \neq 0$):
\begin{align*}
&\left(
\begin{array}{ccc} 
1 & 0 \\
\lambda & 1 \\  \end{array}
\right)s_1^{k'}=
\sum\limits_{\mu \in I_1}{{g}^0_{1,\mu}\otimes \mu^{k'}(1-\lambda\mu)^{r-k'}\xvec}
+ \alpha \otimes (-1)^{k'}\lambda^{r-k'} \yvec ,\\ 
&\left(
\begin{array}{ccc} 
0 & 1 \\
1 & 0 \\  \end{array}
\right)s_1^{k'}=
 \sum\limits_{\mu \in I_1}{{g}^0_{1,\mu} \otimes (-1)^{r+\omega} \mu^{r-k'} \xvec}.
 \end{align*} 
Concluding and expanding the above, a general element of $<s_1^{k'}>_{KZ}$ is given by
\begin{flalign*}
\sum\limits_{\lambda \in I_1}c_{\lambda} &\left(
\begin{array}{ccc} 
1 & 0 \\
\lambda & 1 \\  \end{array}
\right) s_1^{k'} +
c_w \left(
\begin{array}{ccc} 
0 & 1 \\
1 & 0 \\  \end{array}
\right) s_1^{k'}   \\
= &\sum\limits_{\lambda \in I_1} c_{\lambda}\Big( \sum\limits _{s=0}^{r-k'}(-\lambda)^s \binom{r-k'}{s} s_1^{k'+s}
+ (-1)^{k'} \lambda^{r-k'}\alpha \otimes \yvec\Big)
+  
c_w(-1)^{r+\omega}s_1^{r-k'}.
\end{flalign*}
We now bound $\dim<s_1^{k'}>_{KZ}$ in order to use Lemma \ref{lem:FinInd}.
Note that 
by Theorem \ref{thm:lucas}
$\binom{r-p^k(r_k+1)}{s}=\prod\limits_{k,k+1\neq j=0}^{f-1}{\binom{r_j}{s_j}}\binom{p-1}{s_k}\binom{r_{k+1}-1}{s_{k+1}}.$
This implies that the $s$-th summand is zero unless $0 \leq s_j \leq r_j$ for $j\neq k,k+1$, and $0\leq s_k \leq p-1$ and  $0\leq s_{k+1} \leq r_{k+1}-1$.
Additionally, if $0 \leq r_k+1+s_k \leq r_k$ and $k'+s\neq r$ we get that $k'+s < r$ and consequently that $s_1^{k'+s} \in \mathrm{Im}(T)$, allowing us to also ignore elements $s^{k'+s}_1$ where $p-r_k-1 \leq s_k \leq p-1$.
Since $\sum\limits_{\lambda \in I_1}\lambda^s=0$ unless $s=q-1$, taking $c_\lambda = \lambda^{q-1-d}$ where $\binom{r-k'}{d}=0$ or $s_1^{k'+d} \in \mathrm{Im}(T)$ and $c_w=0$ yields 
$q-(p-r_k-1)(r_{k+1})\prod\limits_{k,k+1\neq j=0}^{f-1}(r_j+1)-1$ linear combinations.
Note that for $d=q-1-r+k'$, taking $c_\lambda=(-1)^{k'+1 }\lambda^{d}$ we get that,
\[\sum\limits_{\lambda \in I_1}c_{\lambda} \left(
\begin{array}{ccc} 
1 & 0 \\
\lambda & 1 \\  \end{array}
\right) s_1^{k'} 
=
\sum\limits_{\mu \in I_1}{{g}^0_{1,\mu} \otimes (-\mu)^{r} \xvec}
+ \alpha \otimes \yvec
= T(\mathrm{Id}\otimes \xvec )
\in \mathrm{Im}(T).\]
An additional linear combination is obtained for $c_\lambda=\lambda^{q-1-(r-2k')}$ and $c_w=(-1)^{\omega+1}$.

We get that $\dim <s_1^{p^k(r_k+1)}>_{KZ}= q+1 - \dim \ell_{s_1^{k'}} \leq 
(p-r_k-1)(r_{k+1})\prod\limits_{k,k+1\neq j=0}^{f-1}(r_j+1)$.
By Lemma \ref{lem:char}, we see that the $I$-action on $s_1^{p^k (r_k+1)}$ determines a Serre weight $\tau$ of dimension $(p-r_k-1)(r_{k+1})\prod\limits_{k,k+1\neq j=0}^{f-1}(r_j+1)$.
Using Lemma \ref{lem:FinInd}, we have $<s_1^{p^k (r_k+1)}>_{KZ}\simeq\tau$, and by Proposition \ref{prop:wgtgen1} the statement is true for $n \geq 2$.
If $f=1$, this is exactly \cite[Proposition 2.19]{Sc1}.
\item 
The proof uses the same considerations as in ii.
\end{enumerate}
\end{proof}
\begin{lem}
The image of the $I(1)$-invariant $s_1^r=(-1)^{r+1}\alpha\otimes \yvec$ in $\indkzg\sigma/(T)$  generates a KZ-submodule of length $2^f-1$ and dimension $q+1-\prod\limits_{j=0}^{f-1}(r_j+1)$.
\end{lem}
\begin{proof}
From computations similar to Proposition \ref{prop:invdim}
we conclude that 
\[\dim \ell_{s_1^r}=\dim (\mathrm{Im} (T) \cap <s_1^r>_{KZ}) = \prod\limits_{j=0}^{f-1}(r_j+1).\]
Using Lemma \ref{lem:FinInd} we have $<s_1^r>_{KZ}\simeq \mathrm{Ind}_I^K \chi /\ell_{s_1^r}$, 
implying that $\dim(\mathrm{Ind}_I^K \chi /<s_1^r>_{KZ})$ equals the dimension of the unique irreducible quotient of $\mathrm{Ind}_I^K \chi $, proving the statement.
\end{proof}
\subsection{A basis for the space $(\indkzg\sigma/(T))^{I(1)}$}
We can now give a complete description of $(\indkzg\sigma/T)^{I(1)}$.
For the sets $B_n$ as defined in Section \ref{sec:Not}, set $\qBn^{I(1)}\subset \indkzg\sigma/(T)$  to be the set of elements $\tilde{f}\in \indkzg\sigma/(T)$ which are $I(1)$-invariant and have a lift in $\indkzg\sigma$ which is supported in $B_n$.
\begin{thm} \label{thm:InvBasis}
Recall that $e$ and $f$ are respectively the ramification index and inertia degree of $F/\Qp$, assume $2 < r_j < p-3$ for $0 \leq j \leq f-1$ and define the following sets (for $s_n^k$ and $t_n^k$ as in Definition \ref{def:invs}):

$\mathcal S^l_m=\{s_n^{p^l(r_l+1)}\}_{n\geq m}\bigcup\{ \beta s_n^{p^l(r_l+1)}\}_{n\geq m}$, $~~~\mathcal S_m=\bigcup\limits_{l=0}^{f-1} \mathcal S^l_m$,

$\mathcal T^l_m=\{t_n^{l}\}_{n\geq m}\bigcup\{ \beta t_n^{l}\}_{n\geq m}$, $~~~\mathcal T_m=\bigcup\limits_{l=0}^{f-1} \mathcal T^l_m.$

Then an $I$-eigenbasis for the space $(\indkzg\sigma/(T))^{I(1)}$ of $I(1)$-invariants as an $\fpbar$-vector space is given by the following sets:
\begin{align*}
&\{\mathrm{Id}\otimes x^r , \alpha \otimes y^r\}  & :e=1, f=1 &\text{ ~(Breuil)}\\
&\mathcal S_2\bigcup \{\mathrm{Id}\otimes x^r , \alpha \otimes y^r\}\bigcup \mathcal T_1  & :e>1, f=1 &\text{ ~(Schein)}\\
&\mathcal S_1\bigcup \{\mathrm{Id}\otimes \xvec , \alpha \otimes \yvec\} & :e=1, f>1\\
&\mathcal S_1\bigcup \{\mathrm{Id}\otimes \xvec , \alpha \otimes \yvec\}\bigcup \mathcal T_1& :e>1, f>1
\end{align*}
\end{thm}
\begin{conc} \label{conc:EndScalar}
$\mathrm{End}_G(\indkzg\sigma/(T))\simeq \fpbar$.
\end{conc}
\begin{proof}
The element $\mathrm{Id} \otimes \xvec$ must be sent to a scalar multiple of itself by each endomorphism since by Lemma \ref{lem:char}
there are no other elements in the basis of $(\indkzg\sigma/(T))^{I(1)}$ with the same $I$-eigenvalue.
Since $\mathrm{Id} \otimes \xvec$ generates $\indkzg\sigma/(T)$ as a $G$-module, we are done. 
\end{proof}
The proof of Theorem \ref{thm:InvBasis} is divided into three lemmas. 
Cases $1$ and $2$ were proved by Breuil \cite[Theorem 3.2.4]{B} and Schein \cite[Theorem 2.24]{Sc1} respectively, we prove the new Cases $3$ and $4$. 
\begin{lem} \label{lem:InvBasis1}
Let $\tilde{f}\in B_n\backslash B_{n-1}$ with image in $\overline{B}_n^{I(1)}$, and write $\tilde{f}=\tilde{f}^0_n+\tilde{f}^1_n+f'$ for $\tilde{f}^*_n\in S^k_n$, $*\in \{0,1\}$ and $f'\in B_{n-1}$.
Then we have, 
\[
\tilde{f}^0_n=\sum\limits_{\mu \in I_n}g^0_{n,\mu}\otimes \Big( c(\mu)\xvec +\sum\limits_{k=0}^{f-1}d_k(\mu) \xbutkvec\otimes  x_k^{r_k-1}y_k \Big),\]
\[
\tilde{f}^1_n=\sum\limits_{\mu \in I_n}g^1_{n,\mu}\otimes \Big( c'(\mu)\yvec +\sum\limits_{k=0}^{f-1}d'_k(\mu) \ybutkvec\otimes  y_k^{r_k-1}x_k \Big),
\]
where $c(\mu),c'(\mu),d(\mu),d'(\mu)$
 are specializations at $\mu$ of polynomials $c,c',d,d'\in \fpbar[z_0 , \ldots, z_{n-1}]$ of degree not greater than $q-1$ in each variable $z_j$.
\end{lem}
\begin{proof}
We prove the lemma for $\tilde{f}_n^0$ and $n \geq 0$. 
For $\tilde{f}_n^1$, apply $\beta$ and notice that $\beta \tilde{f}_n^1\in S_n^0$, and $\beta^2$ acts trivially. 

For $n \geq 0$, observe the action of the matrix $\left(
\begin{array}{ccc} 
1 & \pi ^{n} \\
0 & 1 \\  \end{array}
\right)\in I(1)$:

\[\left(
\begin{array}{ccc} 
1 & \pi ^{n} \\
0 & 1 \\  \end{array}
\right)
\left(
\begin{array}{ccc}
\pi^n  & \mu \\
0 & 1 \\  \end{array}
\right)=
\left(
\begin{array}{ccc} 
\pi^n & \mu \\
0 & 1 \\  \end{array}
\right)
\left(
\begin{array}{ccc} 
1  & 1 \\
0 & 1 \\  \end{array}
\right) .\]
Set 
$\tilde{f}^0_n=\sum\limits_{\mu \in I_n}\gn \otimes v_\mu$
where $v_\mu = \sum\limits_{\vec i=\vec 0}^{\vec r}{c_{\vec i}(\mu)\bigotimes\limits_{j=0}^{f-1}x_j^{r_j-i_j}y_j^{i_j}}$.
Using the description of $\mathrm{Im}(T)$, as given in Proposition \ref{prop:inv2}, combined with the fact that $\tilde{f}\in B_n$ and $\tilde{f}_n^0 \in S_n^0$,
we have the following:
\begin{flalign*}
\left(
\begin{array}{ccc} 
1  & \pi ^n \\
0 & 1 \\  \end{array}
\right) \tilde{f}_n^0-\tilde{f}_n^0&
=\sum\limits_{\mu \in I_n}\gn \otimes \left(\left(\begin{array}{ccc} 
1  & 1 \\
0 & 1 \\  \end{array}\right)v_{\mu}-v_{\mu}\right) \in \mathrm{Im}(T) ,
\\
\Longrightarrow A_{\mu} := \left(
\begin{array}{ccc} 
1  & 1 \\
0 & 1 \\  \end{array}
\right) v_{\mu}-v_{\mu}&=\sum\limits_{\vec i=\vec 0}^{\vec r}{c_{\vec i}(\mu)\bigotimes\limits_{j=0}^{f-1}x_j^{r_j-i_j}(x_j+y_j)^{i_j}}-v_\mu \in \fpbar \xvec .
\end{flalign*}
We now show that if a symmetric polynomial of the form 
\begin{equation} \label{Eqn1}
\left(\begin{array}{ccc} 
1  & 1 \\
0 & 1 \\  \end{array}\right)v_{\mu}-v_{\mu}
=
\sum\limits_{\vec i=\vec 0}^{\vec r}{c_{\vec i}(\mu)\bigotimes\limits_{j=0}^{f-1}x_j^{r_j-i_j}(x_j+y_j)^{i_j}}-\sum\limits_{\vec i=\vec 0}^{\vec r}{c_{\vec i}(\mu)\bigotimes\limits_{j=0}^{f-1}x_j^{r_j-i_j}y_j^{i_j}}
\end{equation}
lies in $\fpbar \xvec$, then $c_{\vec i}(\mu)=0$ unless $\vec i= \vec 0$ or $ i_{k}=\delta_{k,k_0}$ for some $k_0$.

Let $\{ {e_k}\}_{k=0}^{f-1}$ be a basis of $\mathbb{R}^f$.
Consider the directed graph $X =(X_v,X_e)$ whose vertices are $X_e=\{\sum\limits_{k=0}^{f-1} a_k e_k : 0 \leq a_k \leq r_k \} \cap \mathbb{Z}^f$ and whose directed edges are $X_v=\{(\vec{i},\vec{j}) : \exists k ,~\vec i = \vec j +  {e_k}\}$.
By a directed path in this graph we mean a sequence of consecutive direct edges.

Every $v\in \mathrm{Sym}^{\vec r}{\fpbar^2}$ can be viewed as a function on $X_e$ taking values in $\fpbar$,
where $\vec i$ corresponds to $\yxvec$. In this interpretation the function corresponding to $A_\mu$ is supported on $\{\vec r\}$.

We see from (\ref{Eqn1}) that a non-zero value of $v_{\mu}$ on a vertex $\vec i$ can only contribute to values of $A_\mu$ on a vertex $\vec j$ if there is a non-trivial directed path from $\vec i$ to $\vec j$, i.e. if and only if $ j_k \leq  i_k$ for all $k$ and $j_{k_0} < i_{k_0}$ for some $k_0$. 

Now, if $\vec i = \vec r- {e_k}$ for some $k$, then the value of $A_\mu$ on $\vec i$ depends only on the value of $v_\mu$ on $\vec r$, but since $A_\mu$ vanishes on every $\vec j \neq \vec r$, we must have $c_{\vec r}(\mu)=0$.
By repeating this argument, we deduce that $v_\mu$ can only take non-zero values on $\vec 0$ and $\{e_k\}_{k=0}^{f-1}$,
i.e. $c_{\vec i}=0$ for $\vec i \notin \{ \vec 0 \} \cup \{e_k\}_{k=0}^{f-1}$,
as these are the only vertices which contribute to nothing besides a value of $A_\mu$ on $\vec r$.

Since $c_{\vec 0}(\mu)$ and $\{c_{e_k}(\mu)\}_{k=0}^{f-1}$ are maps from $I_n$ to $\fpbar$, by 
Lemma \ref{lem:poly} there exist unique polynomials $c$ and $\{d_k\}_{k=0}^{f-1}$ of
degree at most $q-1$ such that $c(\mu)=c_{\vec 0}(\mu)$ and $d_k(\mu)=c_{e_k}(\mu)$ for all $\mu \in I_n$.
\end{proof}

\begin{lem} \label{lem:InvBasis2}
Let $\tilde{f}\in B_n \backslash B_{n-1}$ with image in $\qBn^{I(1)}$, and write $\tilde{f}=\tilde{f}^0_n+\tilde{f}^1_n+f'$ for $\tilde{f}_n^0,\tilde{f}_n^1$ and $f'$ as in Lemma \ref{lem:InvBasis1}.
Then $d_k$ and $d'_k$ are constant for all $0\leq k \leq f-1$.
If $e=1$ then $d_k=d'_k=0$. 
\end{lem}
\begin{proof}
The proof is similar to the case $f=1$ (\cite[Theorem 2.24]{Sc1}) and we prove it by induction. 
Write 
\[
\tilde{f}^0_n=\sum\limits_{\mu \in I_n}g^0_{n,\mu}\otimes \Big( c(\mu)\xvec +\sum\limits_{k=0}^{f-1}d_k(\mu) \xbutkvec\otimes  x_k^{r_k-1}y_k \Big),
\]
and assume inductively that $d_k$ is independent of the parameters $\mu_{n-i}$ for $i\leq m$.
The case $m=0$ is trivial, in order to prove that $d_k$ are independent of $\mu_{n-m-1}$ for $m>0$ we observe that

\[\left(
\begin{array}{ccc} 
1 & \pi ^{n-m} \\
0 & 1 \\  \end{array}
\right)
\left(
\begin{array}{ccc} 
\pi^n  & \mu \\
0 & 1 \\  \end{array}
\right)=
\left(
\begin{array}{ccc} 
\pi^n & \mu' \\
0 & 1 \\  \end{array}
\right)
\left(
\begin{array}{ccc} 
1  & z \\
0 & 1 \\  \end{array}
\right), \]
where $\mu_t = \mu'_t$ for all $0 \leq t< n-m$, with $\mu_{n-m} +1= \mu'_{n-m}$ and $z\in \mathcal O_{F}$. 
Also notice that for $ t > n-m $ the digits of $\mu'$ are different than those of $\mu$, but since we know by induction that $d_k$ is independent of $\mu_{n-t}$ for $t\leq m$, the change of variables from $\mu'$ to $\mu$ does not affect $d_k(\mu)$. 
Now, we can compute: 
\begin{flalign*}
 &\left(
\begin{array}{ccc} 
1 & \pi ^{n-m} \\
0 & 1 \\  \end{array}
\right)\tilde{f}_n^0 - \tilde{f}_n^0 
=\sum\limits_{\mu \in I_n}\gn\otimes \tilde c(\mu)\xvec\\
&+
\sum\limits_{\mu \in I_n}\gn\otimes \sum\limits_{k=0}^{f-1}\Big(d_k([\mu]_{n-m-1},\mu_{n-m-1}-1)-d_k([\mu]_{n-m-1},\mu_{n-m-1})\Big) \xbutkvec\otimes x_k^{r_k-1}y_k . 
\end{flalign*}
Note that $\tilde c(\mu)$ is a polynomial in $\fpbar[\mu_0,\ldots,\mu_{n-1}]$ which depends on $m$, with degrees $(p-1)p^j$ of $\mu_{n-1}$ not appearing for $0 \leq j \leq f-1$.
Using Proposition \ref{prop:inv2} we see that $u^+(\pi^{n-m}) \tilde{f}_n^0 -\tilde{f}_n^0 \in \fpbar\xvec$. 
Since any two elements of the form $\xyvec$ are linearly independent for different indices $i\neq i'$, we must have $d_k(\mu_0,\ldots,\mu_{n-m-1}-1)-d_k(\mu_0,\ldots,\mu_{n-m-1})=0$ for all $k$, and thus $d_k$ is independent of $\mu_{n-m-1}$ and by induction, constant.

Now assume $e=1$, and observe that for $b\in I_1$,
\[\left(
\begin{array}{ccc} 
1 & \pi ^{n-1}[b] \\
0 & 1 \\  \end{array}
\right)
\left(
\begin{array}{ccc} 
\pi^n  & \mu \\
0 & 1 \\  \end{array}
\right)=
\left(
\begin{array}{ccc} 
\pi^n & [\mu]_{n-1}+\pi^{n-1}[\mu_{n-1}+b] \\
0 & 1 \\  \end{array}
\right)
\left(
\begin{array}{ccc} 
1  & P_0(\mu_{n-1},b)\\
0 & 1 \\  \end{array}
\right) ,\]
where $P_0(\mu_{n-1},b)$ is the Witt polynomial 
of Lemma \ref{lem:WittAdd}.
Computation yields,
\begin{flalign} \label{2}
&\left(
\begin{array}{ccc} 
1 & \pi ^{n-1}[b] \\
0 & 1 \\  \end{array}
\right)\tilde{f}_n^0 - \tilde{f}_n^0 
 \\
&=\sum\limits_{\mu \in I_n}\gn\otimes  \Big(c([\mu]_{n-1},[\mu_{n-1}-b])-c([\mu]_{n-1},[\mu_{n-1}])\Big)\xvec 
\notag \\
&+ \sum\limits_{\mu \in I_n}\gn\otimes   \sum\limits_{k=0}^{f-1}d_k\Big( P_0(\mu_{n-1}-b,b)^{p^k} \xvec \Big) \in \mathrm{Im }(T). \notag
\end{flalign}
Recall that for all $0 \leq k \leq f-1$ we have,
\[P_0(\mu_{n-1}-b,b)^{p^k}
=\Big(\sum\limits_{s=1}^{p-1}\mu_{n-1}^{sp^{f-1}}(- b)^{q-sp^{f-1}} \frac{\binom{q}{sp^{f-1}}}{p}\Big)^{p^k}
=\sum\limits_{s=1}^{p-1} \lambda_s(b,k) \mu_{n-1}^{sp^{k-1}},\]
\\where $\lambda_s(b,k) \in \fpbar^\times$.

Assume $d_k$ are non-zero and let 
$
c_{[\mu]_{n-1}}(x)=\sum\limits_{j=0}^{q-1}a_j([\mu]_{n-1})z^{j} \in \fpbar[z]
$  be a polynomial such that 
$c_{[\mu]_{n-1}}(\mu_{n-1})=c(\mu)$, where $a_j \in \fpbar[z_1,\ldots,z_{n-1}]$ are suitable polynomials (which exist by Lemma \ref{lem:poly}). 

Observe that $P_0(b,\mu-b)^{p^k}$ has a non-zero monomial of the form $(-b)^{p^{k-1}}\mu_{n-1}^{(p-1)p^{k-1}}$ 
(if $k = 0$, it has a monomial of the form $(-b)^{p^{f-1}}\mu_{n-1}^{(p-1)p^{f-1}}$).
Write
\begin{flalign} 
&\sum\limits_{j=0}^{q-1} \tilde{a}_j([\mu]_{n-1}) \mu_{n-1}^j
:= 
c_{[\mu]_{n-1}}(\mu_{n-1}-b)-c_{[\mu]_{n-1}}(\mu_{n-1}) \label{3} \\
=&\sum\limits_{j=0}^{q-1} a_j([\mu]_{n-1})\Big(\sum\limits_{s=0}^{j-1}(-b)^{j-s}\binom{j}{s}\mu_{n-1}^{s}\Big).
\label{4}
\end{flalign}
By assumption, (\ref{2}) lies in $\mathrm{Im}(T)$, and since $r_k<(p-1)$, by the description of $\mathrm{Im}(T)$ 
we must have $\tilde{a}_{(p-1)p^{k-1}}([\mu]_{n-1})+(-b)^{p^{k-1}}=0$ for all $b \in \mathbb F_q^{\times}$.
Expanding (\ref{4}), for each $b \in \mathbb F_q^{\times} $ we have an equation,
\[
\tilde{a}_{(p-1)p^{k-1}}([\mu]_{n-1})=\sum\limits_{m=(p-1)p^{k-1}+1}^{q-1} a_m([\mu]_{n-1})(-b)^{m-(p-1)p^{k-1}}\binom{m_{k-1}}{p-1}=-(-b)^{p^{k-1}}.
\]
Multiplying each of these equations by $(-b)^{-p^{k-1}}$, and summing them together,
\begin{equation} \label{5}
\sum\limits_{m=(p-1)p^{k-1}+1}^{q-1} a_m([\mu]_{n-1})\binom{m_{k-1}}{p-1}\sum\limits_{b\in \mathbb F^{\times}_q} (-b)^{m-p^{k}}=\sum\limits_{b\in \mathbb F^{\times}_q}-1.
\end{equation}
The elements on the left hand side of (\ref{5}) are zero except when $m-p^k\equiv  0 \mod q-1$.
In that case, $p^{k-1}<m=p^k$ and so $\binom{m_{k-1}}{p-1}=0$. We conclude the left side of (\ref{5}) equals zero.
Since the right hand side of the equation equals $1$, it must be that $d_k=0$ for all $0 \leq k \leq f-1$. 

The argument holds for $k=0$ as well, recalling that the polynomials we work with are subjected to the relation $x^q=x$.
For $\tilde{f}_n^1$, as in Lemma \ref{lem:InvBasis1}, apply $\beta$ and note that $\beta \tilde{f}_n^1\in S_n^0$ and $\beta^2 \tilde{f}_n^1=\tilde{f}_n^1$.
\end{proof}

\begin{lem} \label{lem:InvBasis3}
Let $\tilde{f}\in B_n \backslash B_{n-1}$ with image in $\qBn^{I(1)}$, and $\tilde{f}_n^0,\tilde{f}_n^1$ and $f'$ as before. 
Then $c(\mu)=\sum\limits_{k=0}^{f-1}\tilde c_k \mu_{n-1}^{p^k(r_k+1)}$ and $c'(\mu)=\sum\limits_{k=0}^{f-1}\tilde c'_k \mu_{n-1}^{p^k(r_k+1)}$ and $\tilde c_k$ and $\tilde c'_k$ are independent of $\mu$, where $c(\mu)$ and $c'(\mu)$ are as defined in Lemma \ref{lem:InvBasis1}.
\end{lem}
\begin{proof}
Assume $f>1$.
Since we know that elements of the form 
$t^s_n=\sum\limits_{\mu \in I_n}{
\gn
\otimes} \xbutsvec \otimes x_s^{r_s-1}y_s$ are $I(1)$-invariant by Proposition \ref{prop:inv1}, we use Lemma \ref{lem:InvBasis2} to assume $\tilde{f}_n^0=\sum\limits_{\mu \in I_n}\gn \otimes c({\mu}) \xvec $.

If $n=0$ the claim follows trivially.
For $n \geq 1$, as in the previous lemma, the following difference lies in $\mathrm{Im}(T)$:
\begin{flalign*}
&\left(
\begin{array}{ccc} 
1 & -\pi ^{n-1} \\
0 & 1 \\  \end{array}
\right)\tilde{f}_n^0 - \tilde{f}_n^0 \\
&=
\sum\limits_{\mu \in I_n}\gn\otimes (\underbrace{c([\mu]_{n-1},[\mu_{n-1}+1])-c([\mu]_{n-1},[\mu_{n-1}])}_\emph{$:=\Delta c$})\xvec \in \mathrm{Im}(T). 
\end{flalign*}

First, note that by the explicit description of Proposition \ref{prop:inv2}, the possible degrees $k$ of $\mu_{n-1}$ that can appear in $\Delta c$  must satisfy $0 \leq k_j \leq r_j$ for all $0 \leq j \leq f-1$.
We claim that the possible degrees $k$, of $\mu_{n-1}$ in $c(\mu)$ must have either $0 \leq k_j \leq r_j$ for all $0 \leq j \leq f-1$, or be of the form $k=p^l(r_l+1)$ for some $0\leq l \leq f-1$.

Otherwise, there is a summand $\mu_{n-1}^k$ in $c(\mu)$ where $k_{j_0}>r_{j_0}$ for some $j_0$
and $k\neq p^{j_0}(r_{j_0}+1)$.
Without loss of generality assume there is no other monomial $\mu_{n-1}^{k'}$ in $c(\mu)$ such that $k_j \leq k'_j $ for all $j$.
We have
\[
(\mu_{n-1}+1)^k-\mu_{n-1}^k
=\sum\limits_{i=0}^{k-1}\prod_{j=0}^{f-1}\binom{k_j}{i_j}\mu_{n-1}^i,
\]
and this implies $\Delta c$ contains all monomials of the form $\mu_{n-1}^{k-p^l}$ where $0 \leq l \leq f-1$ (recall $\binom{k_j}{i_j}=0$ if $i_j > k_j$). In particular $\mu_{n-1}^{k-p^{j_1}}$ appears in $\Delta c$ for $j_1$ such that $j_1\neq j_0$ and $k_{j_1}>0$, and since $k_{j_0} > r_{j_0}$ and $k\neq p^{j_0}(r_{j_0}+1)$, this cannot be.

Thus, using Corollary \ref{cor:lower degrees are trivial}, we can assume $f_n^0$ is of the form 
\[
f_n^0=\sum\limits_{\mu \in I_{n}}g^{0}_{n,\mu}\otimes \left(Q([\mu]_{n-1}) \mu_{n-1}^{r}+\sum\limits_{k=0}^{f-1}\tilde  c_k([\mu]_{n-1}) \mu_{n-1}^{p^k(r_k+1)}\right)\xvec,
\]
 where $\tilde c_k$ and $Q$ depend only on $[\mu]_{n-1}$. 
 Assume $n \geq 2$. By Proposition \ref{prop:inv2} we have $(-1)^r s_1^r+ \alpha \otimes \yvec \in \mathrm{Im}(T)$, so we have  following equality modulo $\mathrm{Im}(T)$ (where we set $\tilde{Q}(\mu)=(-1)^{r+1}Q(\mu)$):
\begin{flalign} 
&\sum\limits_{\mu \in I_{n}}g^{0}_{n,\mu}\otimes Q([\mu]_{n-1})\mu_{n-1}^r\xvec 
=\sum\limits_{\mu \in I_{n-1}}g^{0}_{n-1,\mu}Q(\mu)s_1^r \label{6}
\\
=&\sum\limits_{\mu \in I_{n-1}}g^{0}_{n-1,\mu}\alpha \otimes\tilde{Q}(\mu) \yvec
=\sum\limits_{\mu \in I_{n-1}}g^{0}_{n-2,\mu} \otimes\tilde{Q}(\mu) 
\left(
\begin{array}{ccc} 
1 & [\mu_{n-2}] \\
0 & 1 \\  \end{array}
\right)
\yvec. \label{7}
\end{flalign}
This implies we can replace every term as in the LHS of \ref{6} by a term as in the RHS of \ref{7} in $f$ without changing $\tilde{f}_n^1$. If $n=1$, we can just replace $s_1^r$ by $(-1)^{r+1} \alpha \otimes \yvec$, this does not change $\tilde{f}_1^1$.

We prove inductively similarly to the previous lemma that $\tilde c_k$ are independent of $[\mu]_{n-m}$ for $1 \leq m \leq n-1$. Assume $n \geq 2$, for $n=1$ the statement is clear. 
If $m=1$, this is clear by the definition of $\tilde c_k(\mu)$.
Now, assume that $\tilde c_k$ are independent of $\mu_{n-i}$ for all $i < m$,  and observe that,
\begin{flalign*}
&\left(
\begin{array}{ccc} 
1 & \pi ^{n-m} \\
0 & 1 \\  \end{array}
\right) \tilde{f}_n^0-\tilde{f}_n^0
 \\
&=\sum\limits_{\mu \in I_n}\gn\otimes \sum\limits_{k=0}^{f-1}  \tilde{c}_k ([\mu']_{n-m+1}]){\mu'}_{n-1}^{p^k (r_k+1)} -\tilde{c}_k([\mu]_{n-m+1}) \mu_{n-1}^{p^k (r_k+1)} \xvec,
\end{flalign*}
where $\mu'_j=\mu_j$ for $j< n-m$, with $\mu'_{n-m}=\mu_{n-m}-1$ and for $j>n-m$ we have  $\mu'_j=\mu_j+z_j$  for some $z_j(\mu_{n-m},\dots,\mu_{n-2})\in \fpbar$. 
This is equal to, 
\begin{flalign}
=&\sum\limits_{\mu \in I_n}\gn\otimes \sum\limits_{k=0}^{f-1} \big (\tilde{c}_k ([\mu]_{n-m},\mu_{n-m}-1)-\tilde{c}_k([\mu]_{n-m},\mu_{n-m}) \big)\mu_{n-1}^{p^k (r_k+1)}\xvec \label{8}\\
&+\sum\limits_{\mu \in I_n}\gn\otimes \sum\limits_{k=0}^{f-1} \sum\limits_{j=0}^{r_k} \tilde{c}_k ([\mu]_{n-m},\mu_{n-m}-1)\mu_{n-1}^{j p^k}z_{n-1}^{p^k (r_k+1)- j p^k}\binom{r_k+1}{j}\xvec. \notag
\end{flalign}
Now, $\mu_{n-1}$ can have degree greater than $p^kr_k$ only in the first line of (\ref{8}), and we conclude by Proposition \ref{prop:inv2} that $\tilde{c}_k ([\mu]_{n-m},\mu_{n-m}-1)=\tilde{c}_k ([\mu]_{n-m},\mu_{n-m})$. This implies $\tilde{c}_k$ is independent of $\mu_{n-m}$ as required.

The proof for $\tilde{f}_n^1$ is by applying $\beta$ and redoing this procedure, it does not alter $\tilde{f}_n^0$.

If $f=1$, this was proven in \cite{Sc1}.
\end{proof}

\begin{proof} (of Theorem \ref{thm:InvBasis})
By Lemma \ref{lem:char} we see that the elements of $\mathcal{S}_1$ (and of $\mathcal{T}_1$ if $e>1$) are $I$-eigenvectors, and that those with different $I$-actions must be linearly independent. 
Elements in the two sets above with the same $I$-action are linearly independent because their linear combinations have degrees of the digits $\mu_j$ which cannot appear in $\mathrm{Im}(T)$ by the description given in Proposition \ref{prop:inv2}.

Now, assume $f>1$ and $\tilde{f}\in \indkzg \sigma$ such that its image in the quotient is in $\qBn^{I(1)}$ (where $n\geq 1$ is minimal). 
By the three lemmas $\tilde{f}_n^0$ can be written as a linear combination of elements of the form $s_n^{p^k{(r_k+1)}}$ and $t_n^k$, and the same holds for $\tilde{f}_n^1$ with 
$\beta s_n^{p^k{(r_k+1)}}$ and $\beta t_n^k$.
These elements are $I(1)$-invariant in the quotient, so for $n \geq 2$, the image of $\tilde{f}-\tilde{f}_n^0-\tilde{f}_n^1$ in the quotient is in $\overline{B_{n-1}}^{I(1)}$, and if $n =1 $, then $\tilde{f}-\tilde{f}_1^0-\tilde{f}_1^1$ is supported on $B_0$.
By reiterating this procedure, we arrive at $\tilde{f}_0=\tilde{f}-f'=\mathrm{Id}  \otimes v_0 + \alpha \otimes v_1$ where $f'$ is $I(1)$ invariant in the quotient and supported on $B_n \backslash B_0$, implying that $\tilde{f}_0 \in (\indkzg\sigma / (T))^{I(1)}$.
Now, a simple computation finishes the proof.

\[
\left(
\begin{array}{ccc} 
1 &1 \\
0 & 1 \\  \end{array}
\right)\tilde{f}_0 -\tilde{f}_0 = \mathrm{Id} \otimes \left(
\begin{array}{ccc} 
1 & 1 \\
0 & 1 \\  \end{array}
\right) v_0-v_0 \in \mathrm{Im}(T).
\]
By Proposition \ref{prop:inv2}, this implies the following:
\[
\left(
\begin{array}{ccc} 
1  & 1 \\
0 & 1 \\  \end{array}
\right) v_{0}-v_{0}=\sum\limits_{k=0}^{f-1}c_k (\xbutkvec\otimes x_k^{r_k-1} (x_k + y_k)-\xbutkvec\otimes x_k ^{r_{k-1}}y_k)=0.
 \]
We conclude that $c_k=0$ and $v_0 \in \fpbar \xvec $.
A similar argument holds for $ \alpha \otimes v_1$.
If $f=1$ this was proved in \cite{Sc1}. 
\end{proof}
\subsection{A generalization of the elements $t_n^s$ and their generated Serre weights}
\begin{defn}
Let $\vec k \in \mathbb N_0 ^{f}$ be such that $0 \leq k_j \leq r_j$.
We generalize the elements $t_n^s$ presented before by setting
$t^{\vec k}_n=\sum\limits_{\mu \in I_1}{g_{n,\mu}^0}\otimes \bigotimes \limits _{j=0}^{f-1} x_j^{r_j-k_j}y_j^{k_j}$
and $\widetilde{ \mathcal{T}}^{\vec i}_n=<\{t_n^{\vec i'}\}_{\vec i'\leq\vec i}>_G$.
\end{defn}

Note that $t_1^0=T(\mathrm{Id}\otimes \xvec)$  and $\tilde{ \mathcal{T}}_1^{\vec 0}=T(\indkzg\sigma)$ ($r>0$).

\begin{prop} \label{prop:OtherInvs}
Assume $ 0\leq  k_j \leq \left\lfloor \frac{r_j}{2}  \right\rfloor -1 $ for all $0 \leq j \leq f-1$ and $e>1$. 
Then $t^{\vec k}_n \in \indkzg \sigma / \tilde{ \mathcal{T}}_n^{\vec k-\vec e_m}$ is a non-trivial eigenvector for the $I$-action and in particular $t^{\vec k}_n$ is $I(1)$-invariant for $0\leq m \leq f-1$.
Furthermore, the $KZ$-submodule generated by $t^{\vec k}_n$ is irreducible and isomorphic to the Serre weight $\det^{\omega+k}\otimes \bigotimes \limits_{j=0}^{f-1}\mathrm{Sym}_j^{r_j-2k_j}\fpbar^2$.
\end{prop}
\begin{proof}
We compute the action of $u(a,0)$, $u^+(b)$ and $u^-(\pi c)$ similarly to Proposition \ref{prop:inv1}:
\begin{flalign*}
\left(
\begin{array}{ccc} 
1 & b \\
0 & 1 \\  \end{array}
\right)t_1^{\vec k} -t_1^{\vec k}
&=\sum\limits_{\mu \in I_1}{g}^0_{1,[\mu+b_0]}\otimes \bigotimes\limits_{j=0}^{f-1} x_j^{r_j-k_j}(b_1 ^{p^j}x_j+y_j)^{k_j}-t_1^{\vec k} &\\
&=\sum\limits_{\vec i = \vec 0}^{\vec k}\binom{k}{i} b_1^i \sum\limits_{\mu \in I_1}{g}^0_{1,\mu}\otimes \bigotimes\limits_{j=0}^{f-1} x_j^{r_j-k_j+i_j}y_j^{k_j-i_j}-t_1^{\vec k} \\
&=\sum\limits_{\vec i = \vec 0}^{\vec k}\binom{k}{i} b_1^i t_1^{\vec k - \vec i} - t_1^{\vec k},&\\
\left(
\begin{array}{ccc} 
1 & 0 \\
\pi c & 1 \\  \end{array}
\right)t_1^{\vec k}- t_1^{\vec k}
&=\sum\limits_{\vec i = \vec 0}^{\vec k}\binom{k}{i}  \sum\limits_{\mu \in I_1}{g}^0_{1,\mu}\otimes (-c\mu^2 )^i \bigotimes\limits_{j=0}^{f-1} x_j^{r_j-k_j+i_j}y_j^{k_j-i_j}-t_1^{\vec k},&\\
\left(
\begin{array}{ccc} 
\pi a +1 & 0 \\
0 & 1 \\  \end{array}
\right)t_1^{\vec k} -t_1^{\vec k}
&=\sum\limits_{\vec i = \vec 0}^{\vec k}\binom{k}{i}  \sum\limits_{\mu \in I_1}{g}^0_{1,\mu}\otimes (a \mu )^i \bigotimes\limits_{j=0}^{f-1} x_j^{r_j-k_j+i_j}y_j^{k_j-i_j}-t_1^{\vec k}.&
\end{flalign*}
We continue by induction on $\vec k$.  
For $ \vec k = \vec 0$ we have $t_1^{\vec k}=A_1^0\in (\indkzg\sigma)^{I(1)}$ and it generates a $KZ$-submodules isomorphic to $\sigma$ by Proposition \ref{prop:BasicBasis}. 
The element $t_1^{\vec k}$ is also an eigenvector for the $I$-action.
Now, assume that for all $\vec i < \vec k$ we have $t_1^{\vec i} \in (\indkzg \sigma / \tilde{ \mathcal{T}}_{\vec i-\vec e_m})^{I(1)} $ is an eigenvector for the $I$-action and generates a Serre weight of dimension $\prod\limits_{j=0}^{f-1}({r_j-2i_j+1})$.
In particular, in order to compute $<t_1^{\vec i}>_{KZ}\subset \indkzg \sigma / \tilde{ \mathcal{T}}_{\vec i-\vec e_m}$, we can calculate the possible linear combinations of $t_1^{\vec i}$ over cosets of $K/I$. 
A simple calculation shows that,
\begin{flalign*}
&\left(
\begin{array}{ccc} 
1 & 0 \\
\lambda & 1 \\  \end{array}
\right) t^{\vec i}_1
= \sum\limits_{\nu \in I_1}{{g}^0_{1,\nu} \otimes (1-\lambda \nu)^{r-2i}\otimes \xyvec}+(-1)^i\lambda^{r-2i} \alpha \otimes \yxvec, \\
 &\left(
\begin{array}{ccc} 
0 & 1 \\
1 & 0 \\  \end{array}
\right) t^{\vec i}_1
= \sum\limits_{\nu \in I_1}{{g}^0_{1,\nu} \otimes (-1)^{\omega+r-i}{\nu}^{r-2i}\otimes \xyvec}+(-1)^\omega \alpha \otimes \yxvec.
\end{flalign*}
Observe that $\dim_{\fpbar}<t_1^{\vec i}>_{KZ}=\prod\limits_{j=0}^{f-1}({r_j-2i_j+1})$ which is exactly the number of possible degrees of $\nu$ appearing in the combinations above, taking into account that $\nu^{r-2i} \otimes \xyvec$ and $(-1)^{r-i}\alpha \otimes \yxvec$ always appear with the same coefficient. 
Thus, we must have that $<t_1^{\vec i}>_{KZ}$ is generated by the set ${\mathcal  B}_{t_1^{\vec i}}$ defined as follows (where $\delta_{r-2i,t}$ is the Kronecker delta function),
$$ {\mathcal  B}_{t_1^{\vec i}} =\{\sum\limits_{\nu \in I_1}{g}^0_{1,\nu}\otimes  \nu^{t} \xyvec+(-\delta_{r-2i,t})^{r-i}\alpha \otimes \yxvec\}_{ 0 \leq t_j \leq r_j - 2i_j}.$$

Returning to the actions of $u(a,0)$, $u^+(b)$ and $u^-(\pi c)$, we have that $u^+(b) t_1^{\vec k} - t_1^{\vec k}$ is in $\tilde{ \mathcal{T}}_{\vec k-\vec e_m} $ by assumption, and that $u^-(\pi c) t_1^{\vec k}-t_1^{\vec k}$ is a sum of elements of the form $ \sum\limits_{\nu \in I_1}{g}^0_{1,\nu}\otimes  \nu ^{2i} \bigotimes\limits_{j=0}^{f-1} x_j^{r_j-k_j+i_j}y_j^{k_j-i_j} $. 
Since $2k_j < r_j $ by the conditions of the proposition, we get that $2i_j < r_j -2k_j +2i_j$ and consequently that $u^-(\pi c) t_1^{\vec k}-t_1^{\vec k} \in \tilde{ \mathcal{T}}_{\vec k-\vec e_m} $ since each of the summands presented above is in $<t_1^{\vec k - \vec i}>_{KZ}$ for a suitable $\vec i$. 
As the degrees of $\nu$ in $u(a,0) t_1^{\vec k}-t_1^{\vec k}$ are smaller that those appearing in $u^-(\pi c) t_1^{\vec k}-t_1^{\vec k}$, we conclude that $t_1^{\vec k}$ is an $I(1)$-invariant.
Since $t_1^{\vec k}$ is an $I(1)$-invariant, the $I$-action factors through $I/I(1)=
\bigg\{\left(
\begin{array}{ccc} 
a & 0 \\
0 & d \\  \end{array}
\right)\bigg\}_{a,d \in \mathcal O^{\times}_{F}}
$ and by a simple computation similar to Lemma \ref{lem:char} we get that $i\cdot t_1^{\vec k}= a^{r-2k}(ad)^{k}t_1^{\vec k}$. We denote this character by $\chi_{t_1^{\vec k}}$.
Using Lemma \ref{lem:FinInd}, we see that $\dim_{\fpbar}<t_1^{\vec k}>_{KZ} \geq \dim _{\fpbar} \sigma_{\chi_{t_1^{\vec k}}} $ where $\sigma_{\chi_{t_1^{\vec k}}} $ is the unique Serre weight determined by this character.
By our calculation, $<t_1^{\vec k}>_{KZ}$ is generated by $\mathcal B_{t_1^{\vec k}}$, but since $|\mathcal B_{t_1^{\vec k}}|=\dim _{\fpbar} \sigma_{\chi_{t_1^{\vec k}}}$ using Lemma \ref{lem:FinInd} again we must have $<t_1^{\vec k}>_{KZ}\cong \sigma_{\chi_{t_1^{\vec k}}}$.
For $n \geq 1$, use Proposition \ref{prop:wgtgen1} to conclude that the desired result holds.
\end{proof}

\section{Construction of a universal quotient} \label{sec:4}
\subsection{A universal quotient for an unramified extension of degree 2}
Let $F/\Qp$ be an unramified finite extension with maximal unramified extension $F^{\mathrm{nr}}$. 
Given a continuous, irreducible 
Galois representation $\bar\rho: \mathrm{Gal}(\overline {\mathbb Q}_p / F)\to \mathrm{GL}_2(\fpbar)$ which is sufficiently generic, we associate to it via Serre's weight conjecture the multiset $W(\bar\rho)$ of size $2^f$ of modular Serre weights \cite{BDJ},\cite{GLS}.
Each element in $W(\bar\rho)$ corresponds to a subset $J$ of $\{0,\ldots,f-1\}$, and starting from $\sigma_{\varnothing}\in W(\bar\rho)$ 
we can compute $\sigma_J$. 
An element 
$\sigma=\bigotimes\limits_{\tau\in S}(\det^{w_\tau} \otimes \mathrm{Sym}^{r_{\tau}}\Fq \otimes_{\Fq,\tau}\fpbar)$ where  $S$ is the set of field embeddings $\tau:\mathbb{F}_q \hookrightarrow \fpbar$ 
 is in $W(\bar\rho)$ if and only if $\bar\rho_{|\mathrm{Gal}(\overline{\mathbb{Q}}_p/F^{\mathrm{nr}})}$ is similar to a diagonal matrix $\mathrm{diag}(\varphi,\varphi^q)$,
where $\varphi:\mathrm{Gal}(\overline{\mathbb{Q}}_p/F^{\mathrm{nr}}) \to \fpbar^{\times}$  is a character that can be written as
\[\varphi=\prod\limits_{\tau\in S}\psi_{\tilde \tau}^{(q+1)\omega_{\tau}+r_{\tau}+1},\]
and $\psi_{\tilde \tau}$ is a character of $\mathrm{Gal}(\overline{\mathbb{Q}}_p/F^{\mathrm{nr}})$ obtained through a lift $\tilde\tau$ of $\tau$ to $\mathbb{F}_{q^2}$.
Furthermore, we can take the characters $\psi_{\tilde{\tau}}$ numbered cyclically such that $\psi_{i}^{p}=\psi_{i+1}$.
Note that since we assumed that $2<r_j<p-3$ for all $0 \leq j \leq f-1$, the weights in $W(\bar\rho)$ are determined by the $I$-action on their spaces of $I(1)$-invariants.  
For a more elaborate discussion of this process see \cite[Section 2.4]{Sc4}.

We wish to construct a quotient $U_{\bar\rho}$ of $\indkzg \sigma$, where $\sigma \in W(\bar\rho)$, such that for every supersingular representation $W$ with $\mathrm {soc}_K (W)=\bigoplus\limits_{\sigma \in W(\bar\rho)}\sigma$, and a surjective map $\eta:\indkzg\sigma \twoheadrightarrow W$, the map  $\eta$  factors through $U_{\bar\rho}$.
We first present an example for the case $e=1$ and $f=2$ to illustrate the idea, a general construction is given in Theorem \ref{thm:UnivQuot}.
Throughout this section, we assume $2<r_j<p-3$ for all $0 \leq j \leq f-1$.
\begin{exmpl} \label{ex:QuadExt}
A universal quotient for $e=1$, $f=2$.
\end{exmpl}
We start by calculating the set $W(\bar\rho)$ of modular Serre weights which will appear in the socle of our universal quotient.
\subsubsection*{The set $W(\bar\rho)$ of modular weights for $e=1$ and $f=2$:}
Up to a twist by the determinant, for a generic Galois representation $\bar\rho$ the set $W(\bar\rho)$ is as follows:
\[
W(\bar\rho)= \left\{
\begin{array}{ccc} 
\sigma_{\varnothing} = \mathrm{Sym}_0 ^{r_0}\fpbar ^2 \otimes \mathrm{Sym}_1^{r_1} \fpbar ^2, \\
 \sigma_{\{1\}} =   \mathrm{Sym}_0 ^{r_0-1}\fpbar ^2 \otimes \det^{r_1+1} \otimes \mathrm{Sym}_1^{p-r_1-2} \fpbar ^2, \\
\sigma_{\{0,1\}} = \det^{r_0} \otimes  \mathrm{Sym}_0 ^{p-r_0-1}\fpbar ^2 \otimes \det^{r_1+1} \otimes \mathrm{Sym}_1^{p-r_1-3} \fpbar ^2, \\
 \sigma_{\{0\}} = \det^{r_0} \otimes  \mathrm{Sym}_0 ^{p-r_0-2}\fpbar ^2 \otimes \det^{p-1} \otimes \mathrm{Sym}_1^{r_1+1} \fpbar ^2. 
 \end{array}
\right\}.\] 

We now move to construct a quotient $U_{\bar\rho}$ with $\bigoplus\limits_{\sigma \in W(\bar\rho)}\sigma \subset \mathrm {soc}_{K}(U_{\bar\rho})$.
\subsubsection*{Constructing the universal quotient for $e=1$ and $f=2$:}
Let $V_0 = \indkzg\sigma_\varnothing/(T_{\sigma_{\varnothing}})$.
The element $\mathrm{Id}\otimes \xvec \in V_0$ generates a $KZ$-submodule isomorphic to $\sigma_{\varnothing}$.
We wish to find submodules of $V_0$ which are isomorphic to $\sigma_{\{1\}}$ and $\sigma_{\{0,1\}}$.
By Proposition \ref{prop:invdim} the element $s_1^{p(r_1+1)}$ generates a submodule isomorphic to $\sigma_{\{1\}}$ in $V_0$, so we have a map $\sigma_{\{1\}}\hookrightarrow V_0$ and by Frobenius reciprocity we can obtain a non-zero map $\Phi_{\sigma_{\{1\}}}:\indkzg\sigma_{\{1\}}\to V_0$. 
Now set,
\begin{flalign*}
&s_{\{1\}}=\Phi_{\sigma_{\{1\}}}(\mathrm{Id}\otimes {x_0}^{r_0-1} \otimes {x_1}^{p-r_1-2})
=\sum\limits_{\mu\in I_1} g^0_{1,\mu}\otimes \mu_0^{p(r_1+1)}{x_0}^{r_0} \otimes {x_1}^{r_1}. &\\
&s_{\{0,1\}}=\Phi_{\sigma_{\{1\}}}(\sum\limits_{\mu\in I_1} g^0_{1,\mu}\otimes \mu_0^{r_0} {x_0}^{r_0-1} \otimes {x_1}^{p-r_1-2})
=\sum\limits_{\mu\in I_2} g^0_{2,\mu}\otimes \mu_0^{r_0}\mu_1^{p(r_1+1)}{x_0}^{r_0} \otimes {x_1}^{r_1}.&
\end{flalign*}
We want $s_{\{0,1\}}$ to generate an irreducible $KZ$ module, and thus need $\Phi_{\sigma_{\{1\}}}$ to be compatible with reduction mod $T_{\sigma_{\{1\}}} \in \mathrm{End}_G(\indkzg\sigma_{\{1\}})$. 
Let $V_1=V_0/\Phi_{\sigma_{\{1\}}}(T_{\sigma_{\{1\}}}(\indkzg\sigma_{\{1\}}))$ and consider $\Phi_{\sigma_{\{1\}}}:\indkzg\sigma_{\{1\}}/T_{\sigma_{\{1\}}}\to V_1$. 
Calculating the $I$-action on $s_{\{0,1\}}$, 
\begin{align*}\left(
\begin{array}{ccc} 
a & 0 \\
0 & d \\  \end{array}
\right)s_{\{0,1\}}&=a^r \sum\limits_{\mu\in I_2} g^0_{2,a^{-1}d\mu}\otimes \mu_0^{r_0}\mu_1^{p(r_1+1)}{x_0}^{r_0} \otimes {x_1}^{r_1}&\\
&=a^{p-r_0-1+p(p-r_1-3)}(ad)^{r_0+p(r_1+1)}s_{\{0,1\}}.& 
\end{align*}
Since $\sum\limits_{\mu\in I_1} g^0_{1,\mu}\otimes \mu_0^{r_0} {x_0}^{r_0-1} \otimes {x_1}^{p-r_1-2}$ generates an irreducible $KZ$-submodule in $\indkzg\sigma_{\{1\}}/T_{\sigma_{\{1\}}}$, we have  $<s_{\{0,1\}}>_{KZ}\simeq\sigma_{\{0,1\}}$. 
Now, in a similar way we have $\sigma_{\{0,1\}}\hookrightarrow V_1$ and by Frobenius reciprocity we get a non-zero map $\Phi_{\sigma_{\{0,1\}}}:\indkzg\sigma_{\{0,1\}}\to V_1$. 
We set 
\[s_{\{0\}}=\Phi_{\sigma_{\{0,1\}}}(s_1^{p(p-r_1-2)})=\sum\limits_{\mu\in I_2} g^0_{2,\mu}\otimes \mu_0^{p(p-r_1-2)} \mu_1^{r_0}\mu_2^{p(r_1+1)}{x_0}^{r_0} \otimes {x_1}^{r_1},\]
and calculate the $I$-action:
\begin{align*}
\left(
\begin{array}{ccc} 
a & 0 \\
0 & d \\  \end{array}
\right)s_{\{0\}}&=a^{p-r_0-2+p(r_1+1)}(ad)^{r_0+p(p-1)} s_{\{0\}}.&
\end{align*}
This is the $I$-action on $\sigma_{\{0\}}$, as required.
Now, define $V_2=V_1/\Phi_{\sigma_{\{1\}}}(T_{\sigma_{\{0,1\}}}(\indkzg\sigma_{\{0,1\}}))$, and see that $<s_{\{0\}}>_{KZ}\simeq \sigma_{\{0\}}$ in $V_2$.

For the second phase of our construction take $ U_0 = V_2/(\Phi_{\sigma_{\{0\}}}(T_{\sigma_{\{0\}}}(\indkzg\sigma_{\{0\}})))$, define $M_i=\{\tau\subseteq \mathrm{soc}_K(U_i):\tau \notin W(\bar\rho)\}$ and $N_i=\{\sigma\subseteq \mathrm{soc}_K(U_i):\sigma \in W(\bar\rho)\}$ and set inductively,
\[U_{i+1}=U_i/<\{\Phi_\tau(\indkzg\tau):\tau \in M_i\},\{\Phi_\sigma(T_{\sigma}(\indkzg\sigma)):\sigma \in N_i\}>_G.\]
Finally, take $U_{\bar\rho}=\lim\limits_{\to }U_i$, the universal property presented next will ensure that $U_{\bar\rho}\neq 0$ .
\begin{rem}~
\begin{enumerate}
\item Note that $U_{\bar\rho}$ is still not admissible since we can get infinitely many copies of $\sigma_\varnothing, \sigma_{\{0\}}, \sigma_{\{0,1\}}$ and $\sigma_{\{1\}}$.
\item In the general case the set $N_i$ will be replaced by a set 
which might contain non-modular weights, since in order to obtain all the required weights we will have to pass through weights not appearing in $W(\bar\rho)$.
\end{enumerate}
\end{rem}
\subsubsection*{Proving the universal property of $U_{\bar\rho}$ for $e=1$ and $f=2$:}
Assume there is a surjective map of $G$-modules, $\eta:\indkzg\sigma_{\varnothing} \twoheadrightarrow W$, for $W$ a supersingular representation of $G$ with $\mathrm {soc}_K (W)=\bigoplus\limits_{\sigma \in W(\bar\rho)}\sigma$. 
Note that for an unramified extension of $\Qp$ there exists such $W$ by \cite{BP} 
and in general the construction works but the resulting representations were not proven to be irreducible.
We claim that there exists $\tilde \eta$ such that the following diagram commutes,
\[\begin{tikzpicture}
  \matrix (m) [matrix of math nodes,row sep=3em,column sep=3em,minimum width=2em]
  {
    \indkzg \sigma_{\varnothing} & &W \\
      & U_{\bar\rho} & \\};
  \path[-stealth]
    (m-1-1) 
            edge [-{>[sep= 2pt]>}] node [above] {$\eta$} (m-1-3) 
	 edge  [-{>[sep= 2pt]>}] node [below] { $\varphi$} (m-2-2)
    (m-2-2) 
            edge [->,dashed] node [below] {$\tilde \eta$} (m-1-3);
\end{tikzpicture}\]
Note that $\varphi $ is the reduction map taking an element to its image in the quotient $U_{\bar\rho}$. 
Since $\varphi$ is onto, define $\tilde \eta (x) = \eta (\varphi^{-1}(x))$.
Both $\varphi$ and $\eta$ are $G$-equivariant, so $\tilde \eta$ is well defined if $\ker \varphi \subseteq \ker \eta$. 

We present the proof for $U_0$, showing the universal property for $U_{\bar\rho}$ uses similar arguments and is presented in full generality in the proof of Thereom \ref{thm:UnivQuot}.
Take $\Phi_{\varnothing}$ to be the identity map. 
By the construction of $U_0$, we know that,
\[\ker \varphi_0 = <\{\Phi_{\sigma_{J}} (T_{\sigma_{J}}(\indkzg \sigma_{J}))\}_{J\subseteq \{0,1\}}>_{G},\]
where $\varphi_0:\indkzg\sigma_{\varnothing}\to U_0$ is the quotient map to $U_0$ and $\Phi_{\sigma_{J}} (T_{\sigma_{J}}(\indkzg \sigma_{J}))$ are elements which are lifts of images of $\Phi_{\sigma_{J}} \circ T_{\sigma_{J}}$ in $\indkzg\sigma_{\varnothing}$.

Given a Serre weight $\sigma$, by Proposition \ref{prop:BasicBasis} we know that $\mathrm {soc}_K(\indkzg \sigma)=\bigoplus\limits_{n\in \mathbb N} \sigma $.
Since the operator $T_{\sigma_{\varnothing}}$ is injective, we have  $T_{\sigma_{\varnothing}}(\indkzg \sigma_{{\varnothing}})\simeq \indkzg\sigma_{{\varnothing}}$. 
We must have $(T_{\varnothing}-\lambda_{\varnothing})(\indkzg(\sigma_{\varnothing}))\subseteq \ker \eta$, for some $\lambda_{\varnothing}\in \fpbar$, otherwise, $\mathrm {soc}_K(W)$ will contain two copies of $\sigma_{\varnothing}$ contradicting our assumption.
Since $W$ is supersingular, by the characterization in \cite{BL}, we must have  $\lambda_{\varnothing}=0$, implying that $T_{\sigma_{\varnothing}} (\indkzg \sigma_{\varnothing}) \subseteq \ker \eta$.
Now, since $T_{\sigma_{\varnothing}} (\indkzg \sigma_{\varnothing}) \subseteq \ker \eta$, the map $\Phi_{\sigma_{\{1\}}}$ is defined, and we can use the same argument to conclude that $\Phi_{\sigma_{\{1\}}} (T_{\sigma_{\{1\}}}(\indkzg \sigma_{\{1\}}))\subseteq \ker \eta$. 
Iterating the argument twice more yields $\ker \varphi_0 \subseteq \ker \eta$, implying that $\tilde \eta$ is well defined and the result follows.
In particular, the quotient cannot be zero.
\subsection{A universal quotient for a general unramified extension}
We turn to study the general unramified case. 
We start by establishing an easier way to work with Serre weights.
Recall that for a suitable Galois representation $\bar\rho$, the set of Serre weights $W(\bar\rho)$ corresponds to the power set $P(\{0,\dots,f-1\})$ \cite{Sc4}. 
\begin{lem} 
Set $\sigma_{\varnothing}=\sigma$, for $J \subseteq \{0,\dots,f-1\}$ the parameter and determinant for the Serre weight $\sigma_{J}=\det^{\omega_J} \otimes \bigotimes\limits_{j=0}^{f-1}\mathrm{Sym}_j^{r_j}\fpbar^2$ are given by the function $f_J(r,0)=(r_J,\omega_J)$ with,
\begin{align*}
&r_J= \sum\limits_{j \in J} \big((p-r_j-2) p^j + (-1)^{I_J(j+1)+\delta_{f-1,j}}p^{j+1}\big) +\sum\limits_{j \notin J} r_j p^j,&
\\&\omega_J=\big(\sum\limits_{j\in J} (r_j+1)p^j - p^{j+1}\big)+I_J (f-1) (1-I_J(0)),&
\end{align*} 
where $(r_J,\omega_J) \in (\mathbb {Z} /(q-1) \mathbb {Z})^2 $ and $I_J(x)$ is the indicator function of $J$.
\end{lem}
\begin{proof}
The Serre weight $\sigma_J$ is the weight where 
$\{\psi_{j+I_J (j) f}\}_{0 \leq j \leq f-1}$ is the set of characters appearing in the factorization (\cite{Sc4}, $e=1$):
\begin{equation}\label{factor}
\varphi=\prod\limits_{\tau\in S}\psi_{\tilde \tau}^{r_{\tau}+1}\prod\limits_{\tau\in S}\psi_{\tilde \tau}^{(q_{\mathfrak p}+1)\omega_{\tau}}.
\end{equation}
In order for this to hold, one must start with factorizing $\varphi$ in the following way,
$$\varphi=\prod\limits_{j=0}^{f-1}\psi_{j}^{r_{j}+1}=\prod\limits_{j=0}^{f-1}\psi_{j}^{r_{j}+1}\prod\limits_{j\in J}\psi_{j+f}^{p}\psi_{j+f+1}^{-1}$$
Now, turning to write $\varphi$ in the form of \eqref{factor}, we observe that for $f-1 \neq j \in J$,
 if $j+1\in J$ we have in the $(j+1)$-th place $(\psi_{j+1}\psi_{j+f+1})^{r_{j+1}+1} \psi_{j+f+1}^{p-r_{j+1}-2}$ and if $j+1\notin J$, then $(\psi_{j+1}\psi_{j+f+1})^{-1} \psi_{j+1}^{r_{j+1}+2}$.
If $j= f-1\in J $, we have $0 \in J$ resulting in $(\psi_{0}\psi_{f})^{r_0} \psi_{f}^{p-r_{0}}$ and $0 \notin J$ resulting in $(\psi_{0}\psi_{f})^{0} \psi_{0}^{r_0}$.
Note that since if $f-1 \in J$ and $0 \notin J$ we have $\omega_0=0$, we must add the factor $I_J (f-1) (1-I_J(0))$ to $\omega_J$.
If $j,j+1 \notin J$, we get $(\psi_{j+1}\psi_{j+f+1})^{0} \psi_{j+1}^{r_{j+1}+1}$.
Since $\psi_j^{p}=\psi_{j+1}$ and $\prod\limits_{j=0}^{f-1}(\psi_{j}\psi_{j+f})^{p-1}=1$, we are done since we described explicitly $r_j$ and $\omega_j$ for all $0 \leq j \leq f-1$. 
\end{proof}
\begin{rem} \label{rem:mult}
Note that once $\sigma_{\varnothing}$ is chosen, as we assume $2<r_j<p-3$ for $0\leq j \leq f-1$, we know that $r_J$ determines $\omega_J$ for weights $\sigma_J \in W(\bar \rho)$.
Furthermore, note that for such $r_J$, since $F/\Qp$ is unramified it is evident from the formula that each $\sigma_J \in W(\bar\rho)$ appears with multiplicity one. 
\end{rem}
\begin{defn} \label{def:A_j}
For $0 \leq j \leq f-1$ and $r\in \mathbb Z /((q-1)\mathbb Z)$ define $A_j(r)=r-2 (r_j+1) p^j\mod q-1$.
\\For a set $J\subseteq \{0,\dots,f-1\}$ define ${A_J}_{\downarrow}=\prod\limits_{j\in J} A_j$, where the operators $A_j$ are in descending order (the first operator acting is $A_{\max(J)}$), and ${A_J}_{\uparrow}$ for operators acting in ascending order. Also take $J^c$ to be the complement of $J$ in $\{0,\dots,f-1\}$.
\end{defn}
\begin{lem}\label{lem:Steps1}
Let $J \subseteq \{0,\ldots,f-1\}$ and set $a=\min (J)$ and $J^c_a=\{j\in J^c: j>a\}$. Then: 
\begin{enumerate}
 \item If $0, f-1 \in J$ or $ f-1 \notin J$ then ${A_{J^c}}_{\uparrow} {A_J}_{\downarrow} {A_{J^c}}_{\uparrow} (r_{\varnothing})=r_J.$
 \item If $0 \notin J$ and $f-1 \in J$ then ${A_{J^c_a}}_{\uparrow} {A_J}_{\downarrow} {A_{J^c_a}}_{\uparrow} (r_{\varnothing})=r_J.$
\end{enumerate}
\end{lem}
\begin{proof}
($\textit{1.}$) We apply ${A_{J^c}}_{\uparrow}$ to $r_{\varnothing}$, the operators are taken in ascending order, so if $j\notin J$, we must subtract $p^{j+1}$ if $j+1\in J$ and add $p^{j+1}$ if $j+1\notin J$. 
Also, notice that if $f-1 \notin J$, the action of ${A_{J^c}}_{\uparrow}$ subtracts $p^0$ from the digit in the 0-th place in base $p$. We get,
\[
\sum\limits_{j \notin J} ((p-r_j-2)p^j+(-1)^{I_J(j+1)}p^{j+1})+\sum\limits_{j \in J} r_j p^j-(1-I_J(f-1))((-1)^{I_J(0)}+1).
\]
Now, applying ${A_{J}}_{\downarrow}$, since these operators are applied in descending order, if $f-1\neq j \in J$, we subtract $p^{j+1}$.
Notice that if $f-1,0\in J$, we add $p^0$. The result is,
\begin{align*}
&\sum\limits_{j \notin J} ((p-r_j-2)p^j+(-1)^{2I_J(j+1) }p^{j+1})+\sum\limits_{j \in J} ((p-r_j-2) p^j - p^{j+1})
\\&-(1-I_J(f-1))(1+(-1)^{I_J(0)})+2I_J(f-1)I_J(0).
\end{align*}
Applying ${A_{J^c}}_{\uparrow}$ once more, notice that the term $-(1-I_J(f-1))(1+(-1)^{I_J(0)})$ cancels, arriving at the required expression: 
\begin{align*}
&\sum\limits_{j \notin J} (r_j p^j+(-1)^{I_J(j+1)}p^{j+1}+(-1)^{I_J(j+1)+1}p^{j+1})&
\\&+\sum\limits_{j \in J} ((p-r_j-2) p^j +(-1)^{I_J(j+1)}p^{j+1})+2I_J(f-1)I_J(0)\\
&=\sum\limits_{j \in J} ((p-r_j-2) p^j +(-1)^{I_J(j+1)+\delta_{f-1,j}}p^{j+1})+\sum\limits_{j \notin J} r_j p^j=r_J.&
\end{align*}
The proof of \textit{(2.)} uses similar computations.
\end{proof}
\begin{exmpl} \label{ex:BadWeights}
Unnecessary weights in the construction for $f=3$.
\end{exmpl}
As noted previously, the case $f=2$ is special in the sense that one can construct the desired quotient without passing through weights not in $W(\bar\rho)$, and thus get a quotient whose $K$-socle contains only the irreducible submodules that are in $W(\bar\rho)$.
This breaks for $f\geq 3$, and the weights not in $W(\bar\rho)$ appearing in the construction for $f=3$ and in the $K$-socle of the resulting quotient are given bellow.

In this case, $W(\bar\rho)$ contains $8$ weights, of which $6$ are obtained without passing through weights not in $W(\bar\rho)$;
these are weights corresponding to the subsets $\varnothing,\{2\},\{2,1\},\{2,1,0\},\{1,0\}$ and $\{0\}$.
The undesired weights through which we pass to obtain $\{1\}$ and $\{2,0\}$ are summarized in Table \ref{tab:BadWeights} below (given in the format of parameters $(r_0,\omega_0,r_1,\omega_1,r_2,\omega_2)$).
\begin{table}[h]
\resizebox{12cm}{!}{
\begin{tabular}{ | l | p{11cm} |} 
    \hline
    				& \multicolumn{1}{c|}{Weights not in $W(\bar\rho)$ appearing in the process of obtaining $\sigma_J$}\\ \hline
    $\{1\}$ 		& $(p-r_0-2,r_0+1,r_1-1,0,r_2,0)$,$(p-r_0-3,r_0+1,r_1-1,0,r_2,0)$,$(p-r_0-3,r_0+1,r_1-1,0,p-r_2-2,r_2+1),(p-r_0-3,r_0+1,p-r_1-1,r_1,p-r_2-3,r_2+1),$$(r_0+1,p-1,p-r_1-2,r_1,p-r_2-3,r_2+1)$  \\ \hline
    $\{2,0\}$ 		& $(r_0,0,p-r_1-2,r_1+1,r_2-1,0),(r_0-1,0,p-r_1-2,r_1+1,p-r_2-2,r_2),(p-r_0-1,r_0,p-r_1-3,r_1+1,p-r_2-2,r_2)$ \\ \hline
    \end{tabular}}
\caption[Table caption text]{Unnecessary weights appearing in the $K$-socle of the resulting quotient for $f=3$.}
\label{tab:BadWeights}
\end{table}
\begin{thm} \label{thm:UnivQuot}
Let $F/\Qp$ be a finite unramified extension, $\bar\rho: \mathrm{Gal}(\overline {\mathbb Q}_p / F)\to \mathrm{GL}_2(\fpbar)$ be a Galois representation,  and $\sigma_{\varnothing}\in W(\bar\rho)$, then the following holds:
\begin{itemize}
\item[(1)]  There is an explicit construction of a quotient $U_{\bar\rho}$ of $\indkzg \sigma_{\varnothing}$ such that $\bigoplus\limits_{\sigma \in W(\bar\rho)} \sigma\subseteq \mathrm{soc}_K(U_{\bar\rho})$ and $\mathrm{soc}_K(U_{\bar\rho})$ can only have as irreducible submodules Serre weights with a parameter that appears in a subcomputation of the formulae mentioned in Lemma \ref{lem:Steps1}.
\item[(2)] The quotient $U_{\bar\rho}$ is universal for supersingular representations of $G$ with $ \mathrm {soc}_K (W)=\bigoplus\limits_{\sigma \in W(\bar\rho)}\sigma$; assume there is a surjective map $\eta:\indkzg\sigma \twoheadrightarrow W$, then there exists a map $\tilde \eta$ such that the following diagram commutes,
$$\begin{tikzpicture}
  \matrix (m) [matrix of math nodes,row sep=3em,column sep=3em,minimum width=2em]
  {
    \indkzg \sigma & &W \\
      & U_{\bar\rho} & \\};
  \path[-stealth]
    (m-1-1) 
            edge [-{>[sep= 2pt]>}] node [above] {$\eta$} (m-1-3) 
	 edge  [-{>[sep= 2pt]>}] node [below] { $\varphi$} (m-2-2)
    (m-2-2) 
            edge [->,dashed] node [below] {$\tilde \eta$} (m-1-3);
\end{tikzpicture}$$
where $\varphi $ is the reduction map taking an element in $\indkzg\sigma$ to its projection in $U_{\bar\rho}$. 
\end{itemize}
\end{thm}
\begin{proof}
Set $P_{f}=\{0,\ldots,f-1\}$ and for $J\subseteq P_{f}$ set $a_J$ to be the composition of the operators defined in Lemma \ref{lem:Steps1} 
such that $a_J(r_{\varnothing})=r_J$.
Each such $a_J$ corresponds to a string with letters in $P_{f}$, where the first letter is the operator which acts first on $r_J$.
Also set $L_k=\{j: j \text{ appears in the } k\text{-th place for some } a_J\}$.
We construct the quotient inductively.
For step zero, set $V_0= \indkzg\sigma_{\varnothing}/(T_{\varnothing})$ and let $\varphi_0$ be the projection map from $\indkzg \sigma_{\varnothing}$ to $V_0$.
Since $T_\varnothing$ is injective, we have $\ker\varphi_0=T_{\varnothing}(\indkzg\sigma_{\varnothing})\simeq \indkzg\sigma_{\varnothing}$.
Since $\sigma_{\varnothing}$ appears as a $K$-submodule infinitely many times in $\indkzg\sigma_{\varnothing}$, and $\mathrm{soc}_K(W)$ contains only one copy of $\sigma_{\varnothing}$, we must have  $(T_{\varnothing} - \lambda_{\varnothing})(\indkzg\sigma_{\varnothing}) \subseteq \ker \eta$ for some $\lambda_{\varnothing}\in \fpbar$. 
Since $W$ is supersingular, by the characterization in \cite{BL} we get that $\lambda_{\varnothing}=0$.
This implies $\ker \varphi_0 \subseteq \ker \eta $. 
Note that there exists non-zero supersingular representations $W$  as above by \cite{BP}.

For the first step, by Theorem \ref{thm:InvBasis} for each $j \in L_1$, we have $\sigma_{(J,1)}\simeq <s^{p^j(r_{\varnothing}+1)}>_K\hookrightarrow V_0$ where $\sigma_{(J,1)}$ is the Serre weight with parameter $r_{(J,1)}=A_j(r_\varnothing)$  and $A_j$ the first operator appearing in $a_J$.
By Frobenius reciprocity, we have a map $\Phi_{(J,1)}:\indkzg \sigma_{(J,1)} \to V_0$, and we set $V_1 = V_0/<\{\Phi_{(J,1)}(T_{(J,1)}(\indkzg\sigma_{(J,1)}))\}_{J \subseteq P_{f}}>_G$ and $\varphi_1$ to be the projection map from $\indkzg \sigma_{\varnothing}$ to $V_1$.
Since $\ker \varphi_0 \subseteq \ker \eta$, the map $\eta$ factors through $V_0$, and also that $\mathrm{Im} \Phi_{(1,j)}=<s_1^{p^j({r_{\varnothing}}_j+1)}>_G\subseteq V_0$, where $<s_1^{p^j({r_{\varnothing}}_j+1)}>_G$ contains infinitely many copies of $\sigma_{(J,1)}$ as a $K$-submodule. 
That implies $\Phi_{(J,1)}((T_{(J,1)}-\lambda_{(J,1)})(\indkzg\sigma_{(J,1)})) \subseteq \ker \eta$, for some $\lambda_{(J,1)}\in \fpbar$ since otherwise we have more than one copy of  $\sigma_{(J,1)}$ in $\mathrm{soc}_K(W)$. 
Since $W$ is supersingular, by the characterization in \cite{BL}, we must have $\lambda_{(J,1)}=0$.
We deduce that $\eta$ factors through $V_1$.

Now, define $\sigma_{(J,i+1)}$ to be the Serre weight with parameter $r_{(J,i+1)}=A_j(r_{(J,i)})$ obtained in the $i$-th step for some $j\in L_{i}$ and assume inductively for $1<i<k$ and $j \in L_i$ that there are quotients $V_i = V_{i-1}/<\{\Phi_{(J,i)}(T_{(J,i)}(\indkzg\sigma_{(J,i)}))\}_{J\subseteq P_{f}}>_G$ and also that $\sigma_{(J,i)}\subseteq V_i$ and $\sigma_{(J,i)}\not\subseteq V_s$ for $s<i-1$.
Notice that the last assumption is possible regardless of the induction, since if $\sigma_{(J,i)}\subseteq V_s$ for $s<i$, we can delete the letters in places $s+1$ to $i$ in the string $a_J$ before defining the sets $L_j$ and use $\sigma_{(J,i)}\subseteq V_s$ that was obtained in the $s$-th step.
Also assume that $\eta$ factors through the quotients $V_i$, meaning that $\ker \varphi_i \subseteq \ker \eta$ for projections $\varphi_i:\indkzg \sigma_{\varnothing}\to V_i$.

For $k\leq 2f$ (each $a_J$ has at most $2f$ letters) and $J\subseteq P_{f}$, the maps $\Phi_{(J,k-1)}:\indkzg\sigma_{(J,k-1)}/(T_{J,k-1})\to V_{k-1}$ are defined, and for suitable $j' \in L_{k-1}$ and $j \in L_k$ we have, 
\[\sigma_{(J,k)}=<s_1^{p^j(r_{(k-1,j')}+1)}>_K \hookrightarrow \indkzg\sigma_{(J,k-1)}/(T_{J,k-1})\xrightarrow{\Phi_{(J,k-1)}} V_{k-1}.\]
Note that by Theorem \ref{thm:InvBasis}, $\sigma_{(J,k)}$ is the Serre weight with parameter $A_j(r_{(J,k-1)})$.
Also, we have $\sigma_{(J,k)}\not \subset \ker \Phi_{(J,k-1)}$. 
Otherwise, we would have
\[\sigma_{(J,k)} \subset \Phi_{(J',s)}(T_{(J',s)}(\indkzg\sigma_{(J',s)}))\subseteq V_{s-1}, \]
for some $s < k-1$ and $J'\subseteq P_{f}$, but $V_{k-1}$ is the first quotient in which $\sigma_{(J,k)}$ is obtained.
By Frobenius reciprocity we have a map $\Phi_{(J,k)}:\indkzg\sigma_{(J,k)}\to V_{k-1}$
and define, 
\[V_k = V_{k-1}/<\{\Phi_{(J,k)}(T_{(J,k)}(\indkzg\sigma_{(J,k)}))\}_{J\subseteq P_{f}}>_G.\] 
Notice that $\sigma_{(J,k)}\subseteq V_k$, and also that $\Phi_{(J,k)}(\mathrm{Id} \otimes \bigotimes\limits_{s=0}^{f-1}x_s^{r_{(J,k)}})$ is an $I(1)$-invariant and generates a copy of $\sigma_{(J,k)}$ in $V_{k-1}$. 
By Proposition \ref{prop:wgtgen1}, we then have infinitely many copies of $\sigma_{J,k}$ in $V_{k-1}$, so for some $\lambda_{(J,k)}\in \fpbar$,
\[\Phi_{(J,k)}((T_{(J,k)}-\lambda_{(J,k)})(\indkzg\sigma_{(J,k)})) \subseteq \ker \eta,\] 
since each submodule of $\mathrm{soc}_K(W)$ has multiplicity one by Remark \ref{rem:mult}.
Since $W$ is supersingular, by \cite{BL} we must have  $\lambda_{(J,k)}=0$, meaning that $\eta$ factors through $V_k$.
Note that $\{\sigma_{(J,k)}\}_{J\subseteq P_{f}}$ are the set of Serre weights with parameters $\{ A_{(a_J)_k} \cdot \ldots \cdot A_{(a_J)_1}(r_{\varnothing})\}_{J\subseteq \{0,\ldots,f-1\}}$.
As the strings $a_J$ are of length at most $2f$, after finitely many steps we are done. We denote the resulting quotient by $U_0$.
Note that $V_k\neq 0$ for all $k$ since $0 \neq \sigma_{(J,k)}\subseteq V_k$.

Now, take $S$ to be the set of Serre weights $\{\sigma_{(J,k)}\}$ we obtained in the construction above, and continue with the following process inductively.
Observe $\mathrm{soc}_K(U_i)$, for each $\tau\subseteq \mathrm{soc}_K(U_i)$ by Frobenius reciprocity we have a map $\Phi_{\tau}:\indkzg\tau \to U_i$. 
Set $M_i =\{ \tau\subset \mathrm{soc}_K(U_i) : \tau \notin S\}$, and $N_i =\{ \tau\subset \mathrm{soc}_K(U_i) : \tau \in S$\} and take, \[U_{i+1}=U_{i}/<\{\Phi_{\tau}(\indkzg \tau ): \tau \in M_i\},\{\Phi_{\sigma}(T_{\sigma}(\indkzg \sigma )): \sigma \in N_i\}>_G.\]
Furthermore, since $\tau \not \subseteq \mathrm{soc}_K (W)$ for all $\tau \in M_i$, we must have  $\Phi_{\tau}(\indkzg \tau ) \subseteq \ker \eta$, as the generator of $\tau$ must be sent to zero in $W$. 
Also, as each $\sigma \in N_i$ appears once in $\mathrm{soc}_K (W)$, and since $W$ is supersingular we must have  $\Phi_{\sigma}(T_{\sigma}(\indkzg \sigma )) \subseteq \ker \eta$.
We conclude that $\eta$ factors through $U_i$ for all $i \in \mathbb N_0$.
Since $\{U_i\}_{i \in {\mathbb N}_0}$ forms a directed system, we can take $U_{\bar\rho}=\lim\limits_{\to }U_i= \indkzg \sigma_{\varnothing}/\bigcup\limits_{i=0}^{\infty} \ker \tilde \varphi_i$ where $\tilde \varphi_i :\indkzg\sigma_{\varnothing}\to U_i$ are the projections to $U_i$. 
Furthermore, since $\ker \tilde \varphi_i \subseteq \ker \eta$ for all $i$, the map $\eta$ factors through $U_{\bar\rho}$ by the universal property of the direct limit, and in particular $U_{\bar\rho}\neq 0$.
This finishes our proof.
\end{proof}

\subsection{A universal quotient for $f=2$ and $e<\min\{\frac{r_j}{2}\}$} \label{sec:GenQuot}
Now, set $f=2$ and $e<\min\limits_{0 \leq j \leq f-1}\{\frac{r_j}{2}\}$ (recall we are assuming $2<r_j<p-3$). 
We reproduce the construction presented earlier in this section with a few alterations to obtain a universal quotient $U_{\bar\rho}$ such that every supersingular representation with a prescribed $K$-socle factors through it as before.
Given a continuous, irreducible Galois representation $\bar\rho: \mathrm{Gal}(\overline {\mathbb Q}_p / F)\to \mathrm{GL}_2(\fpbar)$, we can still associate to it a multiset $W(\bar\rho)$, which will be of size $(2e)^f=4 e^2$.
In this case, the elements of $W(\bar\rho)$ are parametrized by pairs $(J,(\delta_0,\delta_1))$ where $J\subseteq \{0,1\}$ and $0 \leq \delta_0,\delta_1\leq e-1$, allowing us to view $W(\bar\rho)$ as a union of $e^2$ lattices. 
As previously, for a thorough discussion of $W(\bar\rho)$ and how to compute it see \cite{BDJ} or \cite[Section 2.4]{Sc4}. 
Note that taking $e<\min\limits_{0 \leq j \leq f-1}\{\frac{r_j}{2}\}$ ensures that  that the conditions for Proposition \ref{prop:OtherInvs} are met and that each weight appears with multiplicity $1$ in $W(\bar\rho)$, allowing us to show that supersingular representations with a suitable $K$-socle indeed factor through our quotient.
\subsubsection*{Calculating the set $W(\bar\rho)$:}
Start with $\sigma_{(\varnothing,(0,0))} = \mathrm{Sym}_0 ^{r_0}\fpbar ^2 \otimes \mathrm{Sym}_1^{r_1} \fpbar ^2\in W(\bar\rho)$. 
The elements of $W(\bar\rho)$ for a fixed $(\delta_0,\delta_1)$ can be computed in a manner similar to as in Example \ref{ex:QuadExt} and should be viewed as perturbations of the weights obtained for that case ($f=2$,$e=1$). The weights in $W(\bar\rho)$ are summarized in Table \ref{tab:weights}:
\begin{table}[h]
\resizebox{14cm}{!}{
\begin{tabular}{ | l | l | p{5cm} |} 
    \hline
    				& \multicolumn{1}{c|}{$(\delta_0,\delta_1)$}\\ \hline
    $\varnothing$ 	& ${\det}^{\delta_0} \otimes\mathrm{Sym}_0 ^{r_0-2\delta_0}\fpbar ^2 \otimes {\det}^{\delta_1} \otimes \mathrm{Sym}_1^{r_1-2\delta_1} \fpbar ^2$ \\ \hline
    $\{1\}$ 		& ${\det}^{\delta_0} \otimes\mathrm{Sym}_0 ^{r_0-2\delta_0-1}\fpbar ^2 \otimes {\det}^{r_1-e+\delta_1+2} \otimes \mathrm{Sym}_1^{p-r_1+2e-2\delta_1-4} \fpbar ^2$  \\ \hline
    $\{0,1\}$ 		& ${\det}^{r_0-e+\delta_0+1} \otimes\mathrm{Sym}_0 ^{p-r_0+2e-2\delta_0-3}\fpbar ^2 \otimes {\det}^{r_1-e+\delta_1+2} \otimes \mathrm{Sym}_1^{p-r_1+2e-2\delta_1-5} \fpbar ^2$ \\ \hline
    $\{0\}$		& ${\det}^{r_0-e+\delta_0+1} \otimes\mathrm{Sym}_0 ^{p-r_0+2e-2\delta_0-4}\fpbar ^2 \otimes {\det}^{p+\delta_1-1} \otimes \mathrm{Sym}_1^{r_1-2\delta_1+1} \fpbar ^2$ \\
    \hline
    \end{tabular}}
\caption[Table caption text]{The weights of $W(\bar\rho)$}
\label{tab:weights}
\end{table}\\
We next move to constructing a suitable quotient containing these weights.
\subsubsection*{Obtaining the required weights and constructing the quotient:}
We will first show that we can obtain a quotient whose $K$-socle contains all the weights in $W(\bar\rho)$ parametrized by $(\varnothing,(\delta_0,\delta_1))$, where $0\leq \delta_0,\delta_1 \leq e-1$. 
From such a quotient we can produce a quotient whose $K$-socle contains all the weights in $W(\bar\rho)$.

As in Proposition \ref{prop:OtherInvs}, set $\widetilde{ \mathcal{T}}^{\vec i}_n=<\{t_n^{\vec i'}\}_{\vec i'\leq\vec i}>_G$, then the element $t_n^{\vec k}$ generates an irreducible $K$-submodule isomorphic to ${\det}^{k_0} \otimes\mathrm{Sym}_0 ^{r_0-2k_0}\fpbar ^2 \otimes {\det}^{k_1} \otimes \mathrm{Sym}_1^{r_1-2k_1} \fpbar ^2$  in $\indkzg \sigma / \tilde{ \mathcal{T}}_{\vec k-\vec e_m}$.
Thus, if we take \[\mathcal{T}=<t_n^{(\delta_0,\delta_1)}: \delta_0+\delta_1=n+1 \text{ and } \max\{\delta_0,\delta_1\}\leq e-1>_G,\] then in $V_{-1}=\indkzg \sigma_{(\varnothing,(0,0))}/<\mathrm{Im}(T_{(\varnothing,(0,0))}),\mathcal{T}>_G$
 the $K$-submodule generated by $t_{\delta_0+\delta_1}^{(\delta_0,\delta_1)}$ is irreducible and isomorphic to $\sigma_{(\varnothing,(\delta_0,\delta_1))}\in W(\bar\rho)$. 

We now wish to follow the first phase of the construction described in Example \ref{ex:QuadExt} for each $(\delta_0,\delta_1)$ to obtain a quotient of $V_0$ containing the desired weights. 
For a given $(\delta_0,\delta_1)$ by Frobenius reciprocity 
we have a non-zero map:
\[\Phi_{(\varnothing,(\delta_0,\delta_1))}:\indkzg \sigma_{(\varnothing,(\delta_0,\delta_1))} \to V_{-1}.\]
Set $V_0=V_{-1}/<\{\mathrm{Im} (\Phi_{(\varnothing,(\delta_0,\delta_1))}\circ T_{(\varnothing,(\delta_0,\delta_1))})\}_{0\leq \delta_0,\delta_1\leq e-1}>_G$,
and in $V_0$ consider the irreducible $K$-submodule generated by $\Phi_{(\varnothing,(\delta_0,\delta_1))}(s_1^{p(r_{\varnothing,(\delta_0,\delta_1)}+1)})$. 
A close inspection shows it generates a $K$-submodule isomorphic to $\sigma_{(\{1\},(\delta_0,e-\delta_1-1))}$, so we again obtain a non-zero map:
\[\Phi_{(\{1\},(\delta_0,e-\delta_1-1))}:\indkzg \sigma_{(\{1\},(\delta_0,e-\delta_1-1))} \to V_{0}.\]
Define $V_1=V_0/<\mathrm{Im}(\Phi_{(\{1\},(\delta_0,e-\delta_1-1))}\circ T_{(\{1\},(\delta_0,e-\delta_1-1))})_{0\leq \delta_0,\delta_1\leq e-1}>_G$,
and notice that in $V_1$ the element $\Phi_{(\{1\},(\delta_0,e-\delta_1-1))}(s_1^{r_{(\{1\},(\delta_0,e-\delta_1-1))}+1})$ generates a $K$-submodule isomorphic to $\sigma_{(\{0,1\},(e-\delta_0-1,e-\delta_1-1))}$.
Iterating this argument again we define $V_2$ in an analogous way, and see that the image of suitable element in $V_2$ generates a $K$-submodule isomorphic to $\sigma_{(\{0\},(e-\delta_0-1,\delta_1))}$. 
Note that we have exactly followed the process of the first phase of the construction in Example \ref{ex:QuadExt}.
Since we have
$W(\bar\rho)=\bigcup\limits_{0 \leq \delta_0,\delta_1 \leq e-1} S_{(\delta_0,\delta_1)}$ where
\[S_{(\delta_0,\delta_1)}=\{\sigma_{(\varnothing,(\delta_0,\delta_1))},\sigma_{(\{0\},(e-\delta_0-1,\delta_1))},\sigma_{(\{0,1\},(e-\delta_0-1,e-\delta_1-1))},\sigma_{(\{1\},(\delta_0,e-\delta_1-1))}\},\]
we see that $\bigoplus\limits_{\sigma \in W(\bar\rho)} \sigma\subseteq \mathrm{soc}_K(V_2)$. 

Finally, we take
\[U_0=V_2/<\mathrm{Im}(\Phi_{(\{0\},(e-\delta_01,\delta_1))}\circ T_{(\{0\},(e-\delta_0-1,\delta_1))})_{0\leq \delta_0,\delta_1\leq e-1}>_G,\] 
and set $M_i=\{\tau\subseteq \mathrm{soc}_K(U_i):\tau \notin W(\bar\rho)\}$ and $N_i=\{\sigma\subseteq \mathrm{soc}_K(U_i):\sigma \in W(\bar\rho)\}$. 
We define inductively quotients,
\[U_{i+1}=U_i/<\{\Phi_\tau(\indkzg\tau):\tau \in M_i\},\{\Phi_\sigma(T_{\sigma}(\indkzg\sigma)):\sigma \in N_i\}>_G,\]
and by taking $U_{\bar\rho}=\lim\limits_{\to }U_i$ we finish the construction.
Once more, since each weight in $W(\bar\rho)$ appears once, if $\indkzg\sigma_{\varnothing,(0,0)} \twoheadrightarrow W$ is a supersingular representation with $\mathrm{soc}_K(W)=W(\bar\rho)$, then $W$ is a quotient of every $U_i$ (from considerations as those presented in Example \ref{ex:QuadExt}) and hence of $U_{\bar\rho}$.
\addtocontents{toc}{\protect\enlargethispage{\baselineskip}}
\bibliographystyle{amsalpha}

\begin{thebibliography}{99}
\addcontentsline{toc}{section}{References}
\bibitem{BS} M. Bardoe and P. Sin, The permutation modules for ${ GL}(n+1,{\mathbb F}_q)$
              acting on ${\mathbb P}^n({\mathbb F}_q)$ and ${\mathbb F}^{n-1}_q$,
\emph{J. London Math. Soc. (2)} Vol 61, no 1, 58-80 (2000).

\bibitem{BLGG} T. Barnet-Lamb, T. Gee and D. Geraghty, Serre weights for rank two unitary groups,
\emph{Math. Ann.} Vol 356 no. 4, 1551-1598 (2013).

\bibitem{BL} L. Barthel and R. Livn\'{e}, Irreducible modular representations of $\mathrm {GL}_2$ of a local field,
\emph{Duke Math. J.} Vol 75, no. 2, 261–292 (1994).

\bibitem{B} C. Breuil, Sur quelques repr\'{e}sentations modulaires et p-adiques de $\mathrm {GL}_2(\Qp)$ I,
\emph{Compositio Math.} Vol 138, 165-188 (2003).

\bibitem{B1} C. Breuil, The emerging $p$-adic Langlands programme,
\emph{Proceedings of the International Congress of Mathematicians} Hyberabad, India (2010).

\bibitem{BP} C. Breuil and V. Paskunas, Towards a modulo p Langlands correspondence for GL2,
\emph{Mem. Amer. Math. Soc.} Vol 216, no. 1016 (2012).

\bibitem{BDJ} K. Buzzard, F. Diamond, and F. Jarvis, On Serre’s conjecture for mod l Galois
representations over totally real fields,
\emph{Duke Math. J.} Vol 155, no. 1, 105-161 (2010).
Proceedings of the ICM, Beijing 2002, vol. 1, 91–97

\bibitem{Fi} N. Fine, Binomial coefficients modulo a prime,
\emph{Amer. Math. Monthly} Vol 54, 589–592 (1947).

\bibitem{GLS} T. Gee, T. Liu and D. Savitt, The weight part of {S}erre's conjecture for {$\mathrm{GL}(2)$},
\emph{Forum Math. Pi}, Vol 3, e2, 52 (2015).


\bibitem{Her} F. Herzig, The classification of irreducible admissible mod p representations of a $p$-adic $\mathrm {GL}_n$, 
\emph{Invent. Math.} Vol 186, no. 2, 373-434 (2011).

\bibitem{Mo2} S. Morra, On some representations of the Iwahori subgroup,
\emph{J. Number Theory} Vol 132, no. 5, 1075-1150 (2012).

\bibitem{O1} R. Ollivier, Le foncteur des invariants sous l'action du pro-$p$-Iwahori de $\mathrm{GL}_2(F)$,
\emph{J. Reine Angew. Math.} Vol 635, 149-185 (2009).


\bibitem{Sc1} M. Schein, On the universal supersingular mod $p$ representations of ${GL2}_(F)$,
\emph{J. Number Theory} Vol. 141, 242-277 (2014).

\bibitem{Sc2} M. Schein, An irreducibility criterion for supersingular mod p representations of $\mathrm {GL}_2(F)$ for totally ramified extensions of $\mathbb Q_p$, 
\emph{Trans. Amer. Math. Soc.} Vol 363, no. 12, 6269-6289 (2011).

\bibitem{Sc3} M. Schein, Weights in Serre's conjecture for Hilbert modular forms: the ramified case,
\emph{ Israel J. Math.} Vol 166, 369-391 (2008).

\bibitem{Sc4} M. Schein, Serre's modularity conjecture,
\emph{Travaux Math.} Vol 23, 139-172 (2013).

\bibitem{Schr} B. Schraen, Sur la pr\'{e}sentation des repr\'{e}sentations supersinguli\'{e}res de $\mathrm {GL}_2(F)$,
\emph{J. Reine Angew. Math.} Vol 704, 187-208 (2015).

\bibitem{SMA} A. Straub, V. H. Moll, T. Amdeberhan, The $p$-adic valuation of k-central binomial coefficients,
\emph{Acta Arith.} Vol 140, no. 1, 31-42 (2009).

\bibitem{Vi1} M.-F. Vign\'{e}ras, Representations modulo $p$ of the $p$-adic group $\mathrm{GL}(2,F)$,
\emph{Compos. Math.} Vol 140 no. 2, 333-358 (2004).

\end{thebibliography}

\end{document}